\documentclass[11pt]{amsart}

\usepackage{subfiles}
\usepackage{euscript}
\usepackage{tabu}
\usepackage[margin=1in]{geometry} 		         

\usepackage{tabularx}	
\usepackage{booktabs}
\usepackage{diagbox}

\usepackage{graphicx}
\usepackage{faktor}
\usepackage{xcolor}
\usepackage{amsmath,amsthm,amssymb}
\usepackage{mathrsfs}
\usepackage{tikz}
\usepackage{tikz-cd}
\usepackage{accents}
\usepackage{upgreek}
\usepackage{enumitem}
\usepackage{bm}
\usepackage{mathtools}
\usepackage[all]{xy}
\usepackage{caption}

\usepackage{url}
\usepackage{float}
\usepackage{todonotes}
\usepackage{bbm}
\usepackage{colonequals}
\usepackage{longtable}

\usepackage[full]{textcomp}
\usepackage[cal=cm]{mathalfa}
\usepackage{xparse}
\usepackage{comment}
\usepackage{cite}
\usepackage{stmaryrd}

\tikzset{
  commutative diagrams/.cd, 
  arrow style=tikz, 
  diagrams={>=stealth}
}

\theoremstyle{definition}
\newtheorem{innercustomthm}{Theorem}
\newenvironment{customthm}[1]
  {\renewcommand\theinnercustomthm{#1}\innercustomthm}
  {\endinnercustomthm}
  
\makeatletter
\def\@tocline#1#2#3#4#5#6#7{\relax
  \ifnum #1>\c@tocdepth 
  \else
    \par \addpenalty\@secpenalty\addvspace{#2}%
    \begingroup \hyphenpenalty\@M
    \@ifempty{#4}{%
      \@tempdima\csname r@tocindent\number#1\endcsname\relax
    }{%
      \@tempdima#4\relax
    }%
    \parindent\z@ \leftskip#3\relax \advance\leftskip\@tempdima\relax
    \rightskip\@pnumwidth plus4em \parfillskip-\@pnumwidth
    #5\leavevmode\hskip-\@tempdima
      \ifcase #1
       \or\or \hskip 1em \or \hskip 2em \else \hskip 3em \fi%
      #6\nobreak\relax
    \dotfill\hbox to\@pnumwidth{\@tocpagenum{#7}}\par
    \nobreak
    \endgroup
  \fi}
\makeatother

\makeatletter
\DeclareRobustCommand{\cev}[1]{%
  \mathpalette\do@cev{#1}%
}
\newcommand{\do@cev}[2]{%
  \fix@cev{#1}{+}%
  \reflectbox{$\m@th#1\vec{\reflectbox{$\fix@cev{#1}{-}\m@th#1#2\fix@cev{#1}{+}$}}$}%
  \fix@cev{#1}{-}%
}
\newcommand{\fix@cev}[2]{%
  \ifx#1\displaystyle
    \mkern#23mu
  \else
    \ifx#1\textstyle
      \mkern#23mu
    \else
      \ifx#1\scriptstyle
        \mkern#22mu
      \else
        \mkern#22mu
      \fi
    \fi
  \fi
}

\makeatother

\usetikzlibrary{calc}
\usetikzlibrary{fadings}
\usetikzlibrary{decorations.pathmorphing}
\usetikzlibrary{decorations.pathreplacing}
\usetikzlibrary{shapes}

\newcounter{marginnote}
\def\marg#1{\stepcounter{marginnote}{\color{violet}\footnotemark[\themarginnote]\footnotetext[\themarginnote]{\bf\color{violet}#1}}}

\DeclareMathAlphabet{\mathpzc}{OT1}{pzc}{m}{it}

\usepackage[pagebackref]{hyperref}
\hypersetup{
  colorlinks   = true,          
  urlcolor     = blue,          
  linkcolor    = purple,          
  citecolor   = blue             
}
\usepackage{cleveref}
\theoremstyle{theorem}
\newtheorem{theorem}{Theorem}[section]
\newtheorem*{theorem*}{Theorem}

\newtheorem{corollary}[theorem]{Corollary}
\newtheorem{lemma}[theorem]{Lemma}
\newtheorem{proposition}[theorem]{Proposition}
\theoremstyle{definition}
\newtheorem{remark}[theorem]{Remark}
\newtheorem*{remark*}{Remark}

\newtheorem*{runningexample*}{Running example}

\newtheorem*{aside*}{Aside}

\newtheorem{definition}[theorem]{Definition}
\newtheorem{example}[theorem]{Example}

\newtheorem{proposition-definition}[theorem]{Proposition-Definition}

\newcommand{\EB}{\mathbf{T}}

\newcommand{\Zcal}{\mathcal{Z}}

\newcommand{\f}{\mathrm{f}}

\newcommand{\Gm}{\mathbb{G}_{\rm{m}}}

\newcommand{\ul}[1]{\underline{#1}}

\newcommand{\bcd}{\begin{center}\begin{tikzcd}}
\newcommand{\ecd}{\end{tikzcd}\end{center}}

\newcommand{\TT}{\mathbf{T}}
\newcommand{\BB}{\mathbf{B}}

\newcommand{\Aaff}{\mathbb{A}}

\newcommand{\PP}{\mathbb{P}}
\newcommand{\OO}{\mathcal{O}}

\newcommand{\ZZ}{\mathbb{Z}}

\newcommand{\HH}{\mathcal{H}}

\newcommand{\Tail}[1]{\mathbf{T}_{#1}}
\newcommand{\Ell}[1]{\mathbf{Ell}_{#1}}

\newcommand{\Lcal}{\mathcal{L}}
\newcommand{\Ocal}{\mathcal{O}}
\newcommand{\Mcal}{\mathcal{M}}
\newcommand{\Hcal}{\mathcal{H}}

\newcommand{\Bcal}{\mathcal{B}}
\newcommand{\Ecal}{\mathcal{E}}

\newcommand{\Ncal}{\mathcal{N}}
\newcommand{\Ccal}{\mathcal{C}}
\newcommand{\Ucal}{\mathcal{U}}

\newcommand{\Pcal}{\mathcal{P}}

\newcommand{\Scal}{\mathcal{S}}

\newcommand{\Cu}{\mathbf{Cu}}

\newcommand{\calO}{\mathcal{O}}

\newcommand{\Spec}{\operatorname{Spec}}

\NewDocumentCommand{\compatibilitydatum}{m m m m m m O{} O{} O{}}{
\begin{equation*} \begin{tikzcd}[ampersand replacement=\&]
  \: \arrow{r} \& {#1} \arrow{r} \arrow{d}{#7} \& {#2} \arrow{r} \arrow{d}{#8} \& {#3} \arrow{r}{[1]} \arrow{d}{#9} \& \: \\
  \: \arrow{r} \& {#4} \arrow{r} \& {#5} \arrow{r} \& {#6} \arrow{r} \& \:
\end{tikzcd} \end{equation*}}

\NewDocumentCommand{\commutingsquare}{m m m m o O{} O{} O{} O{}}{
\begin{equation}\begin{tikzcd}[ampersand replacement=\&] \label{#5}
  #1 \arrow{r}{#6} \arrow{d}{#7} \& #2 \arrow{d}{#8} \\
  #3 \arrow{r}{#9} \& #4
\end{tikzcd}\IfValueTF{#5}{\label{#5}}{} \end{equation}}

\NewDocumentCommand{\cartesiansquare}{m m m m O{} O{} O{} O{}}{
\begin{equation*}\begin{tikzcd}[ampersand replacement=\&]
  #1 \arrow{r}{#5} \arrow{d}{#6} \arrow[dr, phantom, "\square"] \& #2 \arrow{d}{#7} \\
  #3 \arrow{r}{#8} \& #4
\end{tikzcd} \end{equation*}}

\NewDocumentCommand{\cartesiansquarelabel}{m m m m m O{} O{} O{} O{}}{
\begin{tikzcd}[ampersand replacement=\&]
  #1 \arrow{r}{#6} \arrow{d}{#7} \arrow[dr, phantom, "\square"] \& #2 \arrow{d}{#8} \\
  #3 \arrow{r}{#9} \& #4
\end{tikzcd}\IfValueTF{#5}{\label{#5}}{}
}

\NewDocumentCommand{\triangleofspaces}{m m m O{} O{} O{}}{
\begin{tikzcd} [ampersand replacement=\&]
#1 \arrow{r}{#4} \arrow[bend right]{rr}{#5} \& #2 \arrow{r}{#6} \& #3
\end{tikzcd}}

\newcommand{\stProj}{\mathcal{P}\!\operatorname{roj}}
\newcommand{\spec}{\operatorname{Spec}}

\newcommand{\on}{\operatorname}
\newcommand{\oM}{\overline{\mathcal M}}
\newcommand{\tM}{\widetilde{\mathcal M}}

\newcommand{\la}{\lambda}
\newcommand{\lau}{\lambda_{(1)}}
\newcommand{\lad}{\lambda_{(2)}}
\newcommand{\lat}{\lambda_{(3)}}

\newcommand{\de}{\delta}

\newcommand{\Tac}{\mathbf{Tac}}
\newcommand{\Tri}{\mathbf{Tri}}
\newcommand{\Cus}{\mathbf{Cu}}
\newcommand{\Andil}[1]{\marg{(Andrea) #1}}

\newcommand{\subtag}[1]{\tag{\theparentequation.#1}}

\begin{document}
 
\title[Wall-crossing integral Chow rings of $\oM_{1,n\leq 4}$]{Wall-crossing integral Chow rings of $\oM_{1,n\leq 4}$}
\author[L.~Battistella]{Luca Battistella}
\address[Luca Battistella]{Dipartimento di Matematica, Universit\'{a} di Bologna, Italia}
	\email{luca.battistella2@unibo.it}
\author[A.~Di Lorenzo]{Andrea Di Lorenzo}
	\address[Andrea Di Lorenzo]{Dipartimento di matematica, Università di Pisa, Italia}
	\email{andrea.dilorenzo@unipi.it}

\begin{abstract}
We compute the integral Chow rings of $\oM_{1,n}$ for $n=3,4$. The alternative compactifications introduced by Smyth - and studied further by Lekili and Polishchuk - present each of these stacks as a sequence of weighted blow-ups and blow-downs from a weighted projective space. We compute all the integral Chow rings by repeated application of the blow-up formula.



\end{abstract}

\maketitle
\setcounter{tocdepth}{1}

\tableofcontents

\section{Introduction} 
"Moduli spaces are a geometer's obsession" and studying their intersection theory is a most natural problem. They come equipped with a wealth of \emph{tautological classes} which, sometimes but not always, generate the Chow ring (or even cohomology; see \cite{CL} for recent progress and an up-to-date overview of the literature), and the \emph{relations} among which are largely mysterious (but closely linked to enumerative geometry, see for instance \cite{PPZ}). The study of rational Chow rings of moduli spaces of stable curves was initiated by Mumford \cite{Mum}, spurring a beautiful line of research \cite{Fab, Fab2,  Iza, PenevVakil, CL789}. \emph{Integral} Chow rings \cite{EdidinGrahamIntersection,Vistoli-app} are harder to compute, but they encode far more information. 
Keel computed the integral Chow ring of $\oM_{0,n}$ for every $n\geq3$ \cite{Keel}. In recent years, the topic of integral Chow rings of moduli of stable curves has witnessed a growing interest \cite{Lar, dilorenzo2021polarizedtwistedconicsmoduli, DLPV, 
Per}, 
prompting the development of new computational tools, most notably higher Chow groups and patching. In this context, it is natural to extend the stack of stable curves to include some curves with worse singularities and a torus action, which are then excised.

In this project we compute the integral Chow rings of $\oM_{1,n}$ for low values of $n$ by leveraging the geometry of alternative compactifications. Over a field of characteristic $\neq2,3$ we prove:
\begin{customthm}{A}[$\approx$ Corollaries \ref{cor:picture M13} and \ref{cor:picture M14}]
 For $n\leq 5$, the stack $\oM_{1,n}$ can be obtained by a sequence of weighted blow-ups and blow-downs from a weighted projective stack.
\end{customthm}
\begin{customthm}{B}[$\approx$ Theorems \ref{thm:chow M13} and \ref{thm:chow M14}]
For $n=3,4$, the integral Chow ring of $\oM_{1,n}$ is generated by $\lambda$ and boundary divisor classes; we explicitly compute the ideal of relations.
\end{customthm}
The cases $n=5,6$ are analogous, but computationally more challenging. In \cite{BDL2} we carry out the computation of the integral Chow ring of \emph{every} modular compactification of $\Mcal_{1,n},\ n\leq6,$ in the stack of log canonically polarised Gorenstein curves by patching.

\subsection{Wall-crossing integral Chow rings}
The first ingredient of our strategy are Smyth's alternative compactifications of $\Mcal_{1,n}$ \cite{SmythI}: he introduces \emph{Gorenstein} singularities of genus one - of which he produces an explicit list - to replace elliptic bridges with few special points, thus simplifying the boundary of $\oM_{1,n}$ in the spirit of the Hassett--Keel program (log MMP for the moduli space of curves). Typically, elliptic singularities have moduli (\emph{crimping spaces}), resulting in flipping contractions. The maps among these new compact birational models of $\oM_{1,n}$ can be resolved explicitly, 
on one side by adding to the nodal curve the data of a Gorenstein contraction, and on the other by forgetting the stable model \cite{RSPW,SmythIII}:
 \bcd[cramped]
     && \oM_{1,n}(\frac{3}{2})\ar[dr,"\rho_2"]\ar[dl,"\rho_{\frac{3}{2}}"]&& \cdots \ar[dl, "\rho_{\frac{5}{2}}"] \ar[dr, "\rho_{n-2}"] && \\
    \oM_{1,n}\ar[r,"\rho_1"] & \oM_{1,n}(1)\ar[rr,dashed] && \oM_{1,n}(2)\ar[rr,dashed] && \oM_{1,n}(n-2) \ar[r, "\rho_{n-1}"]  & \oM_{1,n}(n-1).
\ecd
The morphism $\rho_{\frac{2m-1}{2}}$ is the ordinary blow-up centred in the locus of elliptic $m$-bridges, whose exceptional divisors are explicit torus bundles parametrising the attaching data for an elliptic $m$-fold point. The exceptional loci of $\rho_m$, instead, can be identified with $\oM_{1,m}(m-1)$.

The second ingredient of our strategy is Lekili's and Polishchuk's identification 
of $\oM_{1,n}(n-1)$ with 
an explicit closed substack of a certain weighted projective stack \cite{LekiliP}. For $n\leq6$, this is the projective stack itself, or can be identified with a Grassmannian, making the Chow ring of $\oM_{1,n}(n-1)$ easy to compute.
Moreover, by applying the criterion for weighted blow-ups of \cite{criterion} we prove the following result, which may be of independent interest.
\begin{theorem*}[\Cref{prop:r dagger}]
If $m\leq 5$, the morphism $\rho_{m}$ is a weighted blow-up of $\oM_{1,n}(m)$ centered at the locus of elliptic $m$-fold points $\mathbf{Ell}_m(m)$, where the weights $(a_0,\ldots,a_m)$ are determined as in \cite[Theorem 1.5.7]{LekiliP} so that $\oM_{1,m}(m-1)\simeq \Pcal(a_0,\ldots,a_m)$.
\end{theorem*}

The third and final ingredient of our strategy is the weighted blow-up formula for integral Chow rings \cite[Corollary 6.5]{ArenaO}: by repeatedly applying it, we climb our way from $\oM_{1,n}(n-1)$ to $\oM_{1,n}$. Under an assumption of surjectivity of the Gysin pullback to the blow-up centre, the Chow ring of the blow-down can be identified with a subring of the Chow ring of the blow-up, which can be computed by elimination theory: this is the most computationally challenging (and least explicit) step in our scheme.

In passing, we compute the integral Chow rings of all $\oM_{1,n}(m)$ for $m<n$, including in particular Schubert's pseudostable curves $\oM_{1,n}^{\textnormal{ps}}=\oM_{1,n}(1)$. See Table \ref{tab:hilb} for their Hilbert series ($n=4$).

\begin{table}[]
    \centering
    \begin{tabular}{c|c|c|c}
         & $\oM_{1,4}(1)$ & $\oM_{1,4}(2)$ & $\oM_{1,4}(3)$ \\ \hline
         $h(t)$ & $1+ 11t + 18t^2 + 11t^3 +t^4$ & $1+ 7t + 7t^2 + 7t^3 +t^4$ & $1+t +t^2 + t^3 +t^4$
    \end{tabular}
    \caption{Hilbert series of Smyth's alternative compactifications of $\Mcal_{1,4}$.}
    \label{tab:hilb}
\end{table}

\subsection{Relation to other work}
 By tensoring with $\mathbb Q$, our results yield presentations for the \emph{rational} Chow rings of the moduli \emph{spaces} $\overline{M}_{1,n}$ of pointed elliptic curves for $n=3,4$. These have been computed by Belorousski in his thesis \cite{Belorousski} (although the presentations therein are not completely explicit). Also, the rational Chow ring can be deduced from combining work of Petersen on the tautological ring in cohomology \cite{Petersen} with the observation that the Chow--K\"unneth generation property (implying cohomology equals Chow) is stable under weighted blow-ups and blow-downs (see \cite{BaeSchmittII} and \cite{CL}). As a sanity-check, we recover the expected Hilbert series:
\[h(A^*(\oM_{1,3}),t)=1+5t+5t^2+1,\qquad h(A^*(\oM_{1,4}),t)=1+ 12t + 23t^2 + 12t^3 +t^4.\]

Our strategy sheds some light on previous work of Inchiostro \cite{Inchiostro}, who presented $\oM_{1,2}$ as a weighted blow-up of a weighted projective space and used this to compute its cohomological invariants. Bishop has computed the integral Chow ring of the interior $\Mcal_{1,n}$ for all $n\leq10$ \cite{BishopM}, and of the Deligne--Mumford compactification for $n=3$ \cite{Bishop}. Newman \cite{Newman} has computed the first higher Chow group of $\Mcal_{1,n}$ for $n\leq 4$.

\subsection{Conventions}\label{sec:convention} We work over a field of characteristic different from $2$ and $3$, because cusps and tacnodes have extra automorphisms otherwise. A curve (resp. an affine curve) over a scheme $S$ is a proper (resp. affine) and flat morphism $C\to S$ whose geometric fibers are connected, reduced, and of dimension one. 

\begin{table}[hbt]
\caption{Table of notations.}
\begin{tabularx}{\textwidth}{@{}p{0.27\textwidth}X@{}}
\toprule
  $C,\omega_C$ & a Gorenstein curve and its dualising line bundle \\
  $\oM_{1,n}(m)$ & Smyth's stack of $m$-stable curves, $0\leq m<n$ (\Cref{sec:smyth}), \newline resp. a blow-up thereof, if  $m\in\frac{1}{2}\ZZ\setminus\ZZ$ (\Cref{def:aligned})\\
  $\widetilde{\Mcal}_{1,m+1}$ & stack of almost $m$-stable curves (\Cref{sec:altCompII})\\
  $S=\{S_1,\ldots,S_\ell\}$ & a set partition of $[n]=\{1,\ldots,n\}$ into $\ell$ parts, of which the last $\ell-k$ parts are singletons \\
  $\Tail{S}(m)$ & locus of elliptic $S$-bridges in $\oM_{1,n}(m)$ (\Cref{def:ellBridge})\\
  $\Ell{S}(m)$ & locus of elliptic $S$-singularities in $\oM_{1,n}(m)$ (\Cref{def:ellSing})\\
  $\Cu_n,\Tac_{S_1|S_2},\Tri_{S_1|S_2|S_3}$ & cuspidal, tacnodal, $3$-elliptic loci with specified markings distribution \\
  $\mathbb{L}_i,\psi_i$ & cotangent line bundle at the $i$-th marking, and its first Chern class\\
  $\Hcal(m),\lambda_{(m)}$ & Hodge line bundle on $\oM_{1,n}(m)$, and its first Chern class\\
  $\Delta_I,\delta_I$ (or $\Delta_{0,I^c},\delta_{0,I^c}$) & divisor of curves with a rational tail marked by $I^c$, and its class \\
  $\tau_{ij}$ & divisor class of curves with a rational tail marked by $\{i,j\}$ in $\oM_{1,4}(2)$ \\
  $\beta_i,\beta_{ij}$ & exceptional divisor classes in $\oM_{1,4}(\frac{3}{2})$ over $\Tac_i$, $\Tac_{ij}$ in $\oM_{1,4}(2)$\\
\bottomrule
\end{tabularx}
\end{table}

\subsection{Acknowledgements} We thank Vladimir Dotsenko, Michele Pernice, Roberto Pirisi, Dhruv Ranganathan, and Angelo Vistoli for helpful conversations. We thank the anonymous referee for their careful suggestions, which have helped us improve the exposition.

\section{Alternative compactifications and weighted blow-ups}

In this section we prove that Smyth's compactifications $\oM_{1,n}(m)$ are all connected by roofs of weighted blow-ups for $m\leq5$ (see \Cref{prop:r dagger} and the following corollaries). In \S \ref{sec:wb} we recall the notion of a weighted blow-up, and a criterion for detecting when a given morphism is a weighted blow-up. In \S \ref{sec:chow wb} we recall the weighted blow-up formula for integral Chow rings. From \S \ref{sec:smyth} on, we describe the birational geometry of $\oM_{1,n}(m)$. We also develop some tools for the computation of fundamental classes of singular loci in the Chow ring.

\subsection{Weighted blow-ups}\label{sec:wb}
Roughly speaking, a weighted blow-up replaces a closed substack $Y\subseteq X$ with a weighted projective bundle over $Y$. With respect to the ordinary case, the weighting offers an extra degree of freedom. Moreover, weighted blow-ups are not schematic. We recall some basic facts, referring the reader to the foundational work of Quek and Rydh \cite{QuekRydh}.

\begin{definition}
 A \emph{weighted embedding} is a filtration $I_0 \supset I_1 \supset \dots \supset I_n \supset \dots$ of ideal sheaves in $\OO_X$  such that:
\begin{itemize}
        \item $I_0 = \OO_X$,
        \item for every $n$, $m \geq 0$ we have $I_n I_m \subset I_{n+m}$, and
        \item locally in the smooth topology on $X$, there exists a sufficiently large positive integer $d$ such that, for all integers $n \ge 1$, $$I_n=\left(I_1^{l_1}I_2^{l_2} \cdots I_d^{l_d} \; : \; l_i \in \mathbb N, \ \sum_{i=1}^d il_i=n\right)$$ in which case we say that $I_\bullet$ \emph{is generated in degrees} $\le d$.
    \end{itemize}
\end{definition}

The Rees algebra $I_\bullet:=\bigoplus I_nt^n \subset \Ocal_X[t]$ is $\mathbb Z$-graded, corresponding to a $\Gm$-action on $\spec_X(I_\bullet)$, which is therefore a cone. Let $I_{+}$ be the ideal generated by $t$, corresponding to the vertex of the cone.
\begin{definition}
    The \emph{weighted blow-up} of $X$ at $I_\bullet$ is 
    $$Bl_{I_\bullet}X:=\stProj_X\left(I_\bullet\right):=[(\spec_X(I_\bullet) \smallsetminus V(I_+))/\Gm] \rightarrow X.$$ 
    Here $\stProj$ is taken in the sense of quotient stacks (as opposed to the usual $\operatorname{Proj}$, that is the schematic quotient). The weighted blow-up \emph{centre} is $Y=V(I_1)\subseteq X$.
\end{definition}
Ordinary blow-ups correspond to taking $I_n=I^n$. In this case, the Rees algebra is generated in degree $1$, and the blow-up is schematic over the base. 

Given a grading of the generators of an ideal $I$, we can define a filtration by defining $I_n\subset I$ to be the elements of degree $\geq n$.
When the filtration $I_\bullet$ comes from a grading, we say that the weighted blow-up is at $Y=V(I_1)$ rather than at $I_\bullet$.

A fundamental example is the weighted blow-up of $\Aaff^m$ at $0$. Given weights $(a_1,\ldots,a_m)$, this defines a filtration of $I=(x_1,\ldots, x_m)$. The weighted blow-up can then be described as 
$$Bl_{0}\Aaff^m \cong \stProj\left(\frac{k[x_1,\dots x_m][s,x_1',\dots x_m']}{(x_1-s^{a_1}x_1',\dots, x_m-s^{a_m}x_m')}\right)$$ with $x_i$ in degree $0$, $x_i'$ in degree $a_i$, and $s$ in degree $-1$. The exceptional divisor in this case is the weighted projective stack
\[\Pcal(a_1,\ldots, a_m) := [\Aaff^m \smallsetminus \{0\} / \Gm ], \]
where $\Gm$ acts on $\Aaff^m$ with weights $(a_1,\ldots, a_m)$. Given a subset $J=\{j_1,\ldots, j_r\} \subset I$ such that $\on{gcd}(a_{j_1},\ldots, a_{j_r})=d$, the vanishing locus of the $J^c$-labelled coordinates 
has generic stabilizer $\mu_d$. The projective stack $\Pcal(a_1,\ldots, a_m)$ is always smooth, whereas its coarse moduli space is in general not smooth as a scheme.

\begin{definition}
 A weighted blow-up $Bl_{Y}X \to X$ is called \emph{regular} when, smooth-locally on $X$, there exists a flat map $X\to \Aaff^m$ such that $f^{-1}(0)=Y$ (in particular, $Y$ is regularly embedded), and such that (again smooth-locally) the weighted blow-up is isomorphic  to $Bl_{0}\Aaff^m \times_{\Aaff^m} X$.
\end{definition}
Given a regular blow-up, the exceptional divisor is the weighted projective bundle over $Y$
\[\Ecal =  [ \Ncal \smallsetminus \Ncal_0 / \Gm ] \]
where $\Ncal$ is the normal bundle of $Y$ (and $\Ncal_0$ its zero section). It inherits a $\Gm$-action that smooth-locally looks like the linear $\Gm$-action on $\Aaff^m$ with weights $(a_1,\ldots, a_m)$.

The following is a generalisation of Artin's criterion for regular weighted blow-up.
\begin{theorem}[\cite{criterion}, Theorem 1.1]\label{criterion}
Let $\widetilde X$ and $Y$ be smooth, tame, and separated Deligne-Mumford stacks, and let $\pi\colon\mathcal E\to Y$ be a weighted projective bundle with positive dimensional fibers. Assume that there is a closed embedding $\mathcal E\hookrightarrow \widetilde X$
as a Cartier divisor, such that $\OO_{\widetilde X}(\mathcal E)_{|\mathcal E}\simeq\OO_{\mathcal E}(-1)\otimes\pi^*L$ for some  line bundle $L$ on $Y$.

Then there is a smooth, tame, and separated Deligne-Mumford stack $X$, with a closed embedding $i\colon Y\hookrightarrow X$ and a morphism $p\colon\widetilde{X}\to X$ that is a regular weighted blow-up with reduced center $Y$. Moreover,
the resulting square is a pushout in algebraic stacks.
\end{theorem}

\subsection{Integral Chow rings and blow-up formulae}\label{sec:chow wb}
The classical blow-up formulae of \cite[\S 6.7]{Fulton} have been generalised to weighted blow-ups by Arena and Obinna \cite{ArenaO}. We will only use a special case, which we recall for the reader's convenience. Let $\widetilde{X}$ be a regular weighted blow-up of a smooth Deligne--Mumford stack $X$ with centre $Y$, where $\iota\colon Y\hookrightarrow X$ denotes the regular embedding.
Suppose moreover that the pullback homomorphism of Chow rings $\iota^*\colon A^*(X) \to A^*(Y)$ is \emph{surjective}. Then we have the following simple formula \cite[Corollary 6.5]{ArenaO}:
\begin{equation}\label{eqn:blow-up formula}
    A^*(\widetilde{X}) \simeq A^*(X)[t]/(t\cdot \ker(\iota^*), Q(t))
\end{equation}
where $t=-[\Ecal]$ is the opposite of the fundamental class of the exceptional divisor, and $Q(t)\in A^*(X)[t]$ is defined by replacing the $t$-constant term of $c_{\textnormal{top}}^{\Gm}(\Ncal)$ with the fundamental class of $Y$. We explain this in more detail.

By assumption, the normal bundle $\Ncal$ of $Y$ in $X$ has a $\Gm$-action, which, smooth-locally on $Y$, is linear. We can therefore consider the $\Gm$-equivariant Chern classes of $\Ncal$, which are elements of the $\Gm$-equivariant Chow ring of $Y$. The latter being endowed with the trivial $\Gm$-action, we have $A^*_{\Gm}(Y) \simeq A^*(Y)[t]$. In particular, we can think of
$c_{\textnormal{top}}^{\Gm}(\Ncal)$
as a polynomial in the equivariant variable $t$. The surjectivity of $\iota^*$ implies that we can lift its coefficients to $X$. In particular, we can find an element $Q'(t)\in A^*(X)[t]$ such that $\iota^*Q'(t) = c_{\textnormal{top}}^{\Gm}(\Ncal)(t)-c_{\textnormal{top}}^{\Gm}(\Ncal)(0)$. Finally, the polynomial $Q(t)$ is defined by adding in the fundamental class of $Y$, i.e. as \[Q(t)=Q'(t) + [Y]\in A^*(X)[t].\]

Putting everything together, we see that, given a presentation for $A^*(X)$, a presentation for $A^*(\widetilde{X})$ can be computed from the following ingredients:
\begin{enumerate}
    \item generators of the ideal $\ker(\iota^*\colon A^*(X) \to A^*(Y))$,
    \item the fundamental class of the center $[Y]\in A^*(X)$, and
    \item the top equivariant Chern class $c_{\textnormal{top}}^{\Gm}(\Ncal_{Y/X})$.
\end{enumerate}

\subsection{Alternative compactifications of $\mathcal{M}_{1,n}$, I}\label{sec:smyth}
Smyth produced some alternative compactifications of $\mathcal{M}_{1,n}$ in \cite{SmythI}. The idea is to trade a higher local/algebraic complexity (worse-than-nodal singularities, such as cusps and tacnodes) for a lower global/combinatorial complexity (subcurves of genus one are required to have many special points; in a smoothing family, those that don't may be contracted to a singularity as above).

We start by recalling the relevant curve singularities. Smyth classified all the isolated Gorenstein singularities of genus one\footnote{The genus $g=\delta-m+1$ of a singularity is an analytic invariant representing the local contribution of the singular point to the arithmetic genus of a projective curve containing it; for instance, an ordinary $m$-fold point, the singularity of the coordinate axes in $\Aaff^m$, is rational.}: there is exactly one analytic type for every number $m$ of branches.
\begin{definition}
   Let $k$ be a field of characteristic different from $2,3$. A $k$-point of a curve $C$ is called an \emph{elliptic} $m$\emph{-fold point} if the analytic germ of $C$ at $p$ is one of the following:
    \[\hat{\calO}_{C,p}\simeq 
    \left\{ \begin{matrix}
       &k \llbracket x,y \rrbracket/(y^2-x^3) & m=1 &\text{ordinary cusp, } A_2 \\
       &k \llbracket x,y \rrbracket/(y^2-yx^2) & m=2 &\text{ordinary tacnode, }A_3 \\
       &k \llbracket x,y \rrbracket/(x^2y-yx^2) & m=3 &\text{planar triple point, }D_4 \\
       &k \llbracket x_1,\ldots,x_{m-1} \rrbracket/I_m & m\geq 4 &\text{$m$-general lines through the origin of }\Aaff^{m-1}
    \end{matrix}\right.
    \]
    where $I_m$ is the ideal generated by the binomials $x_ix_j-x_ix_h$ for all $i,j,h\in[m-1]$.
\end{definition}
Smyth proved that elliptic $\ell$-fold points arise in one-parameter smoothing families exactly when we contract a 
genus one subcurve $E$ of the nodal central fibre $C$, such that 
$\lvert E\cap\overline{C\setminus E}\rvert =\ell$. 
Up to blowing up the markings, then, we can always replace genus one subcurves with $\ell<n$ special points with an elliptic $\ell$-fold singularity. This motivates the following:
\begin{definition}
    Let $0\leq m<n$. An \emph{$m$-stable} $n$-pointed genus one curve over $S$ is \[(\pi\colon C\to S,p_1,\ldots,p_n)\]
    where $C\to S$ is a proper curve (cf. \Cref{sec:convention}) 
    such that:
    \begin{enumerate}
        \item the morphism $\pi\colon C\to S$ is a family of projective curves of arithmetic genus one, whose fibres have only nodes and elliptic $\ell$-fold points as singularities, with $\ell\leq m$;
        \item the morphisms $p_i:S\to C$ are disjoint sections of $\pi$ whose images are contained in its smooth locus; let us denote their union by $\Sigma$, then
        \item\label{cond:level} for every geometric point $s\in S$, and for every connected subcurve of genus one $E\subset C_s$, the latter contains more than $m$ special points (Smyth calls this number the \emph{level} of $E$): 
        $$\on{lev}(E):=|E\cap \overline{(C_s\smallsetminus E)}| + |\Sigma \cap E|>m;$$ 
        \item\label{cond:aut} for every geometric point $s\in S$, we have $H^0(C_s,\Omega_{C_s}^\vee(-\Sigma\cap C_s))=0$.
    \end{enumerate}
\end{definition}
Notice that $0$-stability coincides with the classical Deligne--Mumford stability.

Every Gorenstein curve of arithmetic genus one admits a decomposition into a \emph{minimal elliptic subcurve} (also called the \emph{core}, or \emph{elliptic spine}; this can be an elliptic curve, a nodal wheel of $\PP^1$, or an elliptic $\ell$-fold point with exactly $\ell$ branches $\simeq\PP^1$) and a union of rational trees, attached nodally. The dualising bundle of the core is trivial. Condition \eqref{cond:level} above can be equivalently required of the core only.
Condition \eqref{cond:aut} above implies the finiteness of the automorphism group. It can be phrased combinatorially: every branch of an elliptic $\ell$-fold point must contain at least one marking or node, and at least one branch must contain at least two. Smyth proves the following:
\begin{theorem}[\cite{SmythI, SmythII}]
    Fix $0\leq m<n$. The moduli stack $\oM_{1,n}(m)$ of $n$-pointed, $m$-stable curves of arithmetic genus one is a proper Deligne--Mumford stack over $\Spec(\ZZ[\frac{1}{6}])$; it is smooth if and only if $m\leq5$. Its coarse moduli space is projective, and it arises from the log MMP of $\overline{\mathrm{M}}_{1,n}$.
\end{theorem}
Smoothness follows from the deformation theory of Gorenstein rings in codimension $3$ \cite[Cor. 2.2]{SmythII}. More generally, these stacks are normal, Gorenstein, and smooth in codimension $\leq6$ \cite[\S 1.5]{LekiliP}. So, only in the range $m\leq5$ are the Chow \emph{rings} of Smyth's spaces defined.

\subsection{Alternative compactifications of $\mathcal{M}_{1,n}$, II}\label{sec:altCompII}

Motivated by homological mirror symmetry, Lekili and Polishchuk constructed alternative compactifications of $\mathcal M_{1,n}$ from $A_\infty$-structures. They consider the moduli space $\Ucal^{sns}_{1,m+1}$ of non-special curves $(C,p_1,\ldots,p_{m+1})$, which we will also call \emph{almost-$m$-stable} and denote by $\widetilde{\Mcal}_{1,m+1}$, consistently with \cite{DLPV}. Lekili and Polishchuk compute an explicit normal form for these curves. When restricted to $\Spec(\ZZ[\frac{1}{6}])$, this yields explicit coordinates on Smyth's moduli spaces $\oM_{1,m+1}(m)$.

\begin{definition}[\cite{LekiliP}, Definition 1.4.1]
    An almost $m$-stable $(m+1)$-marked elliptic curve over $S$ consists of the data $(\pi\colon C\to S,\{p_1,\ldots,p_{m+1}\})$ where $\pi$ is a curve of arithmetic genus one, and the $p_i$ are smooth and disjoint sections of $\pi$ such that, for every geometric point $s\in S$:
    \begin{enumerate}
        \item $h^0(C_s,\Ocal_{C_s}(p_i))=1$ for all $i=1,\ldots,m+1$;
        \item and $\Ocal_{C_s}(p_1+\ldots+p_{m+1})$ is an ample line bundle.
    \end{enumerate}
\end{definition}
The total space $\widetilde\Ucal^{sns}_{1,m+1}$ of the $\Gm$-torsor associated to the Hodge line bundle $\widetilde\Hcal=\pi_*\omega_\pi$ over $\widetilde{\Mcal}_{1,m+1}$ is an affine scheme with natural coordinates, by which Lekili and Polishchuk prove:
\begin{proposition}[\cite{LekiliP}, Proposition 1.1.5 and Theorem 1.4.2]\label{prop:presentation Mtilde}
   Let $V_d$ denote the $\Gm$-representation of weight $d$. Then we have the following identifications:
    \begin{enumerate}
        \item [(m=0)] $\widetilde{\mathcal{M}}_{1,1} \simeq [V_4\oplus V_6/\Gm]$;
        \item [(m=1)] $\widetilde{\mathcal{M}}_{1,2} \simeq [V_2\oplus V_3 \oplus V_4/\Gm]$;
        \item [(m=2)] $\widetilde{\mathcal{M}}_{1,3} \simeq [V_1\oplus V_2^{\oplus 2} \oplus V_3/\Gm]$ ;
        \item [(m=3)] $\widetilde{\mathcal{M}}_{1,4} \simeq [V_1^{\oplus 3}\oplus V_2^{\oplus 2}/\Gm]$.
    \end{enumerate}
\end{proposition}
Interestingly, it can be proved a posteriori that the singularities of an almost $m$-stable $m+1$-marked curve are at worst elliptic $(m+1)$-fold points. Indeed, the unique point with $\Gm$-stabiliser in each of the above corresponds to the $(m+1)$-fold point with one marking for every branch, and the stack is nothing but (the quotient by the natural action of the multiplicative group of) the miniversal deformation of such singularity (see also \cite{Stevens}). 

This has two fundamental consequences: first, by removing $0$ from $\widetilde{\Hcal}$, the $\Gm$-quotient is a weighted projective stack; and second, since $0$ represents the only $m$-unstable curve, the quotient can also be identified with Smyth's $\oM_{1,m+1}(m)$.

\begin{proposition}[\cite{LekiliP}, Corollary 1.5.5 and Theorem 1.5.7]\label{prop:presentation M}
    We have the following isomorphisms:
    \begin{enumerate}
        \item[(m=0)]  $\overline{\mathcal{M}}_{1,1} \simeq \Pcal(4,6)$;
        \item [(m=1)] $\overline{\mathcal{M}}_{1,2}(1) \simeq \Pcal(2,3,4)$;
        \item [(m=2)] $\overline{\mathcal{M}}_{1,3}(2) \simeq \Pcal(1,2,2,3)$ ;
        \item [(m=3)] $\overline{\mathcal{M}}_{1,4}(3) \simeq \Pcal(1,1,1,2,2)$.
    \end{enumerate}
    Moreover, these isomorphisms identify the Hodge line bundle $\Hcal(m)$ with $\Ocal_{\Pcal}(1)$.
\end{proposition}

\subsection{...and maps between them}\label{sec:resolution}
Since the locus of smooth pointed elliptic curves is a dense open substack of every $\oM_{1,n}(m)$, these spaces are all birational. The birational maps:
\begin{equation}\label{eq:comparison_map}
    \oM_{1,n}(m-1)\dashrightarrow\oM_{1,n}(m)
\end{equation}
are regular if and only if $m=1$ (contraction of elliptic tails to ordinary cusps, cf. \cite{Schubert,HH09}) or $m=n-1$ \cite[Corollary 4.5]{SmythII}.

A simultaneous resolution of indeterminacy, inspired by tropical geometry, has been devised by Ranganathan, Santos-Parker, and Wise \cite{RSPW}.

For the purpose of our computations, it is more convenient to have a step-by-step resolution, performing only the necessary blow-ups every time. Such a resolution has been described by Smyth in \cite[\S 2.2]{SmythIII}. The map \eqref{eq:comparison_map} above is not defined precisely where the minimal subcurve of genus one has exactly level $m$. In order to distinguish the components of these loci, it is useful to refine the notion of level of a subcurve to take values in (set-)partitions of $\{1,\ldots, n\}$. It was the intuition of \cite{BKN} that this is the way to obtain all the compactifications of $\mathcal{M}_{1,n}$ by Gorenstein curves with smooth, distinct markings.
Set $[n]:=\{1,\ldots, n\}$.
\begin{definition}
Let $(C,p_1,\ldots,p_n)$ be a Gorenstein curve of genus one.
 \begin{enumerate}
  \item Let $E$ be a subcurve of $C$ of genus one. The level of $E$ is the partition of $[n]$ obtained by grouping all the markings that belong to the same connected component of $C\setminus E$; by convention, markings lying on $E$ make up singletons of $\on{lev}(E)$.
  \item If $C$ contains an elliptic $m$-fold point $q$, the level of $q$ is the partition of $[n]$ describing how the markings subdivide among the connected components of $C\setminus\{q\}$.
 \end{enumerate}
\end{definition}

Let $S\in {[n]\choose m}$ be a partition of $n$ into $m$ parts $S_1,\ldots,S_m$. Up to reordering, suppose that $|S_i|=1$ if and only if $i>k=k(S)$. We adopt the notation $S_1|\ldots|S_m$ to represent $S$.
\begin{definition}\label{def:ellBridge}
     We denote by $\EB_S(m-1)$ the locus of elliptic $m$-bridges of level $S$ in $\oM_{1,n}(m-1)$, i.e. the locus of curves formed by nodally attaching an $(m-1)$-stable elliptic spine marked by $\cup_{i>k} S_i$, with $k$ stable rational curves $C_1,\ldots,C_k$, the subcurve $C_i$ being marked by $S_i$.
\end{definition}
In particular, we have an isomorphism:
\[\EB_S(m-1)=\oM_{1,[k]\sqcup S_{k+1}\sqcup\ldots\sqcup S_m}(m-1)\times\prod_{i=1}^k\oM_{0,\{\star\}\sqcup S_i}\simeq\oM_{1,m}(m-1)\times\prod_{i=1}^k\oM_{0,\lvert S_i\rvert+1}.\]
\begin{definition}\label{def:ellSing}
    Let $\mathbf{Ell}_S(m)$ be the locus of elliptic $m$-fold points of level $S$ in $\oM_{1,n}(m)$. 
\end{definition}

Since this is the worst allowed singularity in this moduli stack, the universal family of curves over $\mathbf{Ell}_S(m)$ is equisingular, in particular all of these curves admit a common normalisation, which is a (disconnected) pointed rational curve. We thus have a morphism:
\[\mathbf{Ell}_S(m)\longrightarrow \oM_{0,S}=\prod_{i=1}^k\oM_{0,\{\star_i\}\sqcup S_i}.\]
What is missing in order to reconstruct the elliptic $m$-fold point is the datum of its local ring as a subalgebra of the semilocal ring of the normalisation. The parameter space for these data is called the \emph{crimping space} in \cite{vdW}, and the moduli space of \emph{attaching data} in \cite[\S 2.2]{SmythII}. Smyth \cite[Lemma 2.2]{SmythI} showed that, for an elliptic $m$-fold singularity, the embedding of $\hat\OO_{C,q}$ in $\oplus_{i=1}^m\hat{\OO}_{C^\nu,\star_i}$ is determined by the choice of a hyperplane in $\oplus_{i=1}^m \mathfrak{m}_{\star_i}/\mathfrak{m}_{\star_i}^2$ not containing any of the lines $\mathfrak{m}_{\star_{i}}/\mathfrak{m}^2_{\star_i}$; this hyperplane is identified with the pullback of the cotangent space to the elliptic singularity under the normalisation map.
Dually, the hyperplane corresponds to a line in $\oplus_{i=1}^m (\mathfrak{m}_{\star_i}
/\mathfrak{m}_{\star_i}^2)^\vee$ not contained in any of the hyperplanes $\on{Ann}(\mathfrak{m}_{\star_i}
/\mathfrak{m}_{\star_i}^2)$. The elliptic singularity can be obtained from its seminormalisation $C^{sn}$ (an ordinary, rational $m$-fold point) by pinching along this line.
Smyth \cite[\S 2.2]{SmythII} explains how to \emph{modularly} compactify this algebraic torus bundle into a projective bundle by \emph{sprouting} the universal curve along the boundary.
\begin{proposition}[{\cite[Prop. 2.16]{SmythII}}]
    There is an isomorphism $\mathbf{Ell}_S(m)\simeq \PP_{\oM_{0,S}}\left(\oplus_{i=1}^k\mathbb L_{\star_i}^\vee\right)$, where $\mathbb L_{\star_i}$ denotes the (pullback of the) cotangent line bundle at the marking $\star_i$ on the $i$-th factor of the product.
\end{proposition}

The rational morphism $\oM_{1,n}(m-1)\dashrightarrow\oM_{1,n}(m)$ flips the locus $\EB_S(m-1)$ to the locus $\mathbf{Ell}_S(m)$.
\begin{definition}\label{def:aligned}
 Let $\oM_{1,n}(m-\frac{1}{2})$ denote the ordinary blow-up of $\oM_{1,n}(m-1)$ in $\EB_m(m-1)=\sqcup_{S\in{[n]\choose m}}\EB_S(m-1)$, and let $\mathbf{Exc}(m-\frac{1}{2})$ denote the exceptional divisor.
\end{definition}

Over $\oM_{1,n}(m-\frac{1}{2})$, Smyth \cite[Proposition 2.3]{SmythIII} constructs a diagram of curves:
\bcd
&\Ccal(m-\frac{1}{2})\ar[dl,"\sigma"]\ar[dr,"\tau"]&\\
\Ccal(m-1)&&\Ccal(m)^\prime
\ecd
where $\sigma$ is the blow-up of the sections $p_{S_{k+1}},\ldots,p_{S_m}$ over $\mathbf{Exc}(m-\frac{1}{2})$, and $\tau$ is the contraction of the strict transform of the core, i.e. 
the minimal subcurve of arithmetic genus one, over $\mathbf{Exc}(m-\frac{1}{2})$. The resulting curve $\Ccal(m)^\prime$ is $m$-stable, as $(m-1)$-stability implies that at least one part $S_i$ has cardinality at least two, hence it induces a morphism:
\[\rho_{m}\colon \oM_{1,n}\left(m-\frac{1}{2}\right)\to \oM_{1,n}(m).\]
We state our first result. We denote by $\mathbf{Ell}_m(m)$ the locus of curves in $\oM_{1,n}(m)$ with the worst possible singularity, that is an elliptic $m$-fold point; $\mathbf{Ell}_m(m)$ is the union of all the irreducible substacks $\Ell{S}(m)\subset \oM_{1,n}(m)$ with $|S|=m$.
\begin{theorem}\label{prop:r dagger}
If $m\leq 5$, the morphism $\rho_{m}$ is a weighted blow-up of $\oM_{1,n}(m)$, centered at $\mathbf{Ell}_m(m)$, with weights $(a_0,\ldots,a_m)$, where $\oM_{1,m}(m-1)\simeq \Pcal(a_0,\ldots,a_m)$.
\end{theorem}
We postpone the proof of this theorem to \S \ref{subsec:proof}. Next, we describe the first three special cases. The result below recovers Inchiostro's description of $\oM_{1,2}$ \cite{Inchiostro}.
\begin{corollary}
The moduli stack of stable $2$-pointed curves of genus one $\oM_{1,2}$ is the weighted blow-up of $\oM_{1,2}(1)\simeq\Pcal(2,3,4)$ with weights $(4,6)$ in the point $\mathbf{Ell}_{[2]}(1)=:\mathbf{Cu}_2$ representing the $2$-pointed cusp.
\end{corollary}

\begin{corollary}\label{cor:picture M13}
    There is a diagram:
    \bcd
    \oM_{1,3}\ar[r,"\rho_1"] & \oM_{1,3}(1)\ar[r,"\rho_2"] & \oM_{1,3}(2)\simeq\Pcal(1,2,2,3),
    \ecd
    where $\rho_2$ is a blow-up with weights $(2,3,4)$ in the three points $\mathbf{Ell}_{i,j|k}(2)=:\mathbf{Tac}_k$ parametrising tacnodal curves, and $\rho_1$ is the blow-up with weights $(4,6)$ in the curve $\mathbf{Ell}_{[3]}(1)=:\mathbf{Cu}_3\simeq\oM_{0,4}$ parametrising cuspidal curves.
\end{corollary}

\begin{corollary}\label{cor:picture M14}
    There is a diagram:
    \bcd
    &&  \oM_{1,4}(\frac{3}{2})\ar[dl,"\rho_{\frac{3}{2}}"]\ar[dr,"\rho_{2}"]&& \\
    \oM_{1,4}\ar[r,"\rho_1"] & \oM_{1,4}(1)\ar[rr,dashed] && \oM_{1,4}(2)\ar[r,"\rho_3"] & \oM_{1,4}(3)\simeq\Pcal(1,1,1,2,2),
    \ecd
    where 
    \begin{itemize}[leftmargin=.5cm]
        \item $\rho_3$ is a blow-up with weights $(1,2,2,3)$ in the six points $\mathbf{Ell}_{h,i|j|k}(3)=:\mathbf{Tri}_{h,i}$ parametrising curves with planar triple points;
        \item $\rho_2$ is the blow-up with weights $(2,3,4)$ in the four curves $\mathbf{Ell}_{h|i,j,k}(2)=:\mathbf{Tac}_h\simeq\oM_{0,4}$ parametrising tacnodal curves with a component marked only by $p_h$, and in the three curves $\mathbf{Ell}_{h,i|j,k}(2)=:\mathbf{Tac}_{h,i}\simeq\PP^1$ parametrising tacnodal curves with a component marked only by $p_h$ and $p_i$;
        \item $\rho_{\frac{3}{2}}$ is the ordinary blow-up in the three surfaces $\EB_{h,i|j,k}(1)\simeq \oM_{1,2}(1)$ of elliptic bridges;
        \item $\rho_1$ is the blow-up with weights $(4,6)$ in the surface $\mathbf{Ell}_{[4]}(1)=:\mathbf{Cu}_4\simeq\oM_{0,5}$ of cuspidal curves.
    \end{itemize}
\end{corollary}

\subsection{Universal line bundles and tautological classes}
From a flat family of projective, Gorenstein curves with smooth, distinct markings, we can define the following bundles over the base:
\[ \mathbb L_i(\pi\colon C\to S, p_1,\ldots,p_n) := p_i^*\omega_{\pi},\quad \Hcal(\pi\colon C\to S, p_1,\ldots,p_n):=\pi_*\omega_{\pi}.\]
The first one is the cotangent line bundle at the $i$-th marking, and the second one is called the Hodge bundle. Their first Chern classes will be denoted by $\psi_{i,\pi}$, respectively by $\lambda_\pi$. A morphism between moduli spaces is often induced by a morphism between their universal curves, in which case we will be able to compare the bundles above; to avoid confusion we will decorate them as appropriate: for instance, we shall use $\Hcal(m)$ and $\lambda_{(m)}$ for the universal $m$-stable curve over $\oM_{1,n}(m)$ (unless the number of markings plays a role).

Let us now specialise to the case of arithmetic genus one. In this case, the Hodge bundle has rank one, and the evaluation map $\Hcal=\pi_*\omega_\pi\to p_i^*\omega_\pi=\mathbb L_i$ is generically an isomorphism. Indeed, we can be more precise: by pushing down the short exact sequence
\[0\to \omega_\pi(-p_i)\to\omega_\pi\to\omega_{\pi|p_i}\to 0,\]
and noticing that $\on{R}^1\pi_*\omega_\pi(-p_i)\to\on{R}^1\pi_*\omega_\pi$ satisfies cohomology and base-change, we conclude that the evaluation map fails to be surjective precisely when $h^0(\OO(p_i))=h^1(\omega(-p_i))>h^1(\omega)=1$, that is, when $p_i$ does not lie on the minimal subcurve of genus one (it is special, in the terminology of Lekili and Polishchuk).

\begin{proposition}[{\cite[Lemma 1.1.1]{LekiliP}}]\label{prop:psiVSlambda}
    Let $(\pi\colon C\to S, p_1,\ldots,p_n)$ be a pointed Gorenstein curve of genus one. Let $\Delta_{0,i}$ denote the locus in $S$ where the $i$-th marking $p_i$ lies on a rational tail, and let us assume that $\Delta_{0,i}$ is a Cartier divisor in $S$, whose class we denote by $\delta_{0,i}$. Then:
    \begin{equation}
        \psi_{i,\pi}=\lambda_{i,\pi}+\delta_{0,i}.
    \end{equation}
    In particular, over $\oM_{1,m+1}(m)$, every $\mathbb L_i$ is identified with $\Hcal(m)$.\footnote{This generalises the fact that elliptic curves are parallelisable.}
\end{proposition}

Let us now assume that we have two families $\pi_i\colon C_i \to S$ of pointed Gorenstein curves of genus one over a base $S$, related by a birational contraction $\rho\colon C_1\to C_2$ - an isomorphism generically on $C_2$, and with connected fibres. Assume that the exceptional locus $E$ of $\rho$ is a vertical Cartier divisor, i.e. a subcurve of $C_1$ over a Cartier divisor $\Delta$ in the base (whose class we denote by $\delta$) such that $C_1$ is regular around $\Sigma:=E\cap\overline{C_{1|\Delta}\setminus E}$. We shall compare the Hodge line bundles associated to $C_1$ and $C_2$; in order to do this, it is enough to assume that $S$ is the spectrum of a discrete valuation ring. Denote by $\omega_i$ the relative dualising sheaf of $\pi_i$, and let $\mathcal{H}_i=\pi_{i*}\omega_i$ be the Hodge line bundle associated to $C_i$, whose class we denote by $\lambda_i$. We recall Smyth's definition of a \emph{balanced} subcurve of genus one, borrowing the notation $E$ and $\Sigma$ from above: the condition is that the chains of rational curves within $E$ separating the nodes of $\Sigma$ from the core of $E$ all have the same length (for a more general, not necessarily regular, one-parameter family, being balanced can be phrased in terms of the tropical distance from the core to the components adjacent $E$). 
\begin{proposition}\label{prop:lambdaVSlambda}
The following hold true:
    \begin{enumerate}
        \item If $\Ecal$ is a rational tail, then $\rho^*\omega_2=\omega_1(-\Ecal)$. If $\Ecal$ is a rational bridge, then $\rho^*\omega_2=\omega_1$. In both cases, the $\lambda$-classes agree: $\lambda_1=\lambda_2$.
        \item If $\Ecal$ is a genus one subcurve, then it is balanced, and $\rho^*\omega_2=\omega_1(\Ecal)$. So, $\lambda_2=\lambda_1+\delta$.
    \end{enumerate}
\end{proposition}
\begin{proof}
    The first fact goes back to Knudsen's thesis \cite[Lemma 1.6]{KnudsenII}. The second one is basically the first claim of \cite[Prop. 3.6 and Cor. 3.7]{SmythII}: consider the short exact sequence
\[0\to \omega_{1}\to\rho^*\omega_2\to\rho^*\omega_{2|\Ecal}=\OO_{\Ecal}\to 0;\]
pushing it down along $\pi_1=\pi_2\circ\rho$, and applying the projection formula (since $\rho$ is a birational contraction, we have $\rho_*\OO_{C_1}=\OO_{C_2}$), we obtain:
\[0\to\Hcal_1\to\Hcal_2\to\OO_S(\Delta)\to 0,\]
where right-exactness follows because the $\on{R}^1\pi_{1,*}\omega_{C_i}$, for $i=1,2$, are isomorphic (by Serre duality, both of them can be identified with $\OO_S$).
\end{proof}

We will also use the following fact \cite[Proposition 3.4]{SmythII}, comparing the $\lambda$-class of a one-parameter family of $\ell$-elliptic singularities with the $\psi$-classes of its pointed normalisation:
\begin{proposition}\label{prop:SmythLambda}
    Let $B\to \oM_{1,n}(m)$ be a one-parameter family corresponding to a curve $C\to B$ with a section $q\colon B\to C$ of elliptic $\ell$-fold points (same $\ell$ for all $b\in B$). Let $(\widetilde C,q_1,\ldots,q_\ell)\to B$ be the pointed 
    normalisation, where $q_1,\ldots, q_\ell$ are the preimages of $q$. Then $\la_{(m)}(C).B=-\psi(\widetilde C,q_i).B$ for any $i\in[\ell].$
\end{proposition}

\subsection{Proof of \Cref{prop:r dagger}}\label{subsec:proof}
The locus of elliptic $m$-bridges $\EB_m(m-1)$ decomposes into the disjoint union of $\EB_S(m-1)$, each of which is a closed regular substack of codimension $k=k(S)$, with split normal bundle:
\begin{equation}
    N_{\EB_S/\oM_{1,n}(m-1)}=\bigoplus_{i=1}^k\mathbb L^\vee_{1,i}\boxtimes \mathbb L^\vee_{0,\star_i}.
\end{equation}
Now, recall from Proposition \ref{prop:psiVSlambda} that on the first factor, which is isomorphic to $\oM_{1,m}(m-1)$, we may identify every cotangent line bundle $\mathbb L_{1,i}$ with the Hodge line bundle $\Hcal(m-1)$, so
\begin{equation}
    N_{\EB_S/\oM_{1,n}(m-1)}=\Hcal(m-1)^\vee\boxtimes\left(\bigoplus_{i=1}^k \mathbb L^\vee_{0,\star_i}\right).
\end{equation}

Recall that $\Ell{m}(m) \subset \oM_{1,n}(m)$ parametrises curves with an elliptic $m$-fold point: this locus decomposes into the disjoint union of the components $\Ell{S}(m)$ with $\#S=m$. We may similarly decompose the exceptional divisor of the blow-up $\rho_{m}\colon \oM_{1,n}\left(m-\frac{1}{2}\right)\to \oM_{1,n}(m)$ into a disjoint union indexed by partitions, and we obtain the following identification:

    \begin{equation}\label{eq:exceptional}
    \mathbf{Exc}_S\left(m-\frac{1}{2}\right)\simeq \oM_{1,m}(m-1)\times\PP_{\oM_{0,S}}\left(\bigoplus_{i=1}^k \mathbb L^\vee_{0,\star_i}\right),
    \end{equation}
    with normal bundle
    \[\mathcal N_{\mathbf{Exc}(m-\frac{1}{2})/\oM_{1,n}(m-\frac{1}{2})}=\mathbb \Hcal(m-1)^\vee \boxtimes \OO(-1).
    \]
    The morphism $\rho_{m}$, restricted to $\mathbf{Exc}_S(m-\frac{1}{2})$, coincides with the projection onto the second factor of \eqref{eq:exceptional}.
    Thanks to Proposition \ref{prop:presentation M}, we can conclude by applying the criterion for smooth weighted blow-downs \cite[Theorem 1.1]{criterion}. Finally, observe that given any weighted blow-up $Bl_Y X \to X$, the exceptional divisor is (locally on the center $Y$ of the weighted blow-up) isomorphic to a product of the base with a weighted projective stack $\mathcal{P}(a_0,\ldots, a_m)$. The weights of the blow-up are precisely the values $a_0,\ldots, a_m$. They can be recovered from the orders of the automorphism groups of the torus-fixed points. Therefore, we can deduce the weights of $\rho_m$ by looking at the weights appearing in the description of $\oM_{1,m}(m-1)$ as a weighted projective stack. This concludes the proof of \Cref{prop:r dagger}.

\subsection{Fundamental classes and normal bundles}\label{sec:tilda}
From \Cref{sec:chow wb} and \Cref{sec:resolution} we know that we need to compute the fundamental class and the normal bundle of the elliptic $m$-fold loci $\Ell{S}(m)$ in $\oM_{1,n}(m)$.
It is easy to perform such a calculation in $\widetilde{\Mcal}_{1,m}$. We would like to pull it back along the map forgetting $n-m$ markings. In general this is only a rational map, but for $n\leq 4$ it is actually regular. We defer a general discussion of forgetful maps to \Cref{appendix}.


We start from Lekili and Polishchuk's description of the space of almost $m$-stable curves. 

\begin{lemma}\label{lm:class of Z}
    For $m\leq 3$, the Chow ring of $\widetilde{\Mcal}_{1,m}$ is isomorphic to $\ZZ[t]$, where $t$ is the first Chern class of $\widetilde{\Hcal}(m)$. Moreover, the point $\Zcal$ representing the elliptic $m$-fold point has class:
    \begin{enumerate}
        \item $\Ncal_{\Zcal/\widetilde{\Mcal}_{1,m}}=\widetilde{\Hcal}(m)^{\otimes4}\oplus\widetilde{\Hcal}(m)^{\otimes6}$ and $[\Zcal] = 24t^2$ when $m=1$;
        \item $\Ncal_{\Zcal/\widetilde{\Mcal}_{1,m}}=\widetilde{\Hcal}(m)^{\otimes2}\oplus\widetilde{\Hcal}(m)^{\otimes3}\oplus\widetilde{\Hcal}(m)^{\otimes4}$ and $[\Zcal] = 24t^3$ when $m=2$;
        \item $\Ncal_{\Zcal/\widetilde{\Mcal}_{1,m}}=\widetilde{\Hcal}(m)^{\otimes1}\oplus\widetilde{\Hcal}(m)^{\otimes2}\oplus\widetilde{\Hcal}(m)^{\otimes2}\oplus\widetilde{\Hcal}(m)^{\otimes3}$ and $[\Zcal] = 12t^4$ when $m=3$.
    \end{enumerate}
\end{lemma}
\begin{proof}
    From the descriptions of \Cref{prop:presentation Mtilde}, we deduce that $\widetilde{\Mcal}_{1,m}$ is a vector bundle over $\Bcal\Gm$. Therefore, by homotopy invariance of Chow groups we get
    \[A^*(\widetilde{\Mcal}_{1,m+1}) \simeq A^*(\Bcal\Gm) \simeq \ZZ[t].\]
    We also have $[\Zcal]=[\{0\}]_{\Gm}$, i.e. the fundamental class of the zero section of a vector bundle over $\Bcal\Gm$. The latter is equal to the top Chern class of the vector bundle itself. Using the fact that $c_i^{\Gm}(V_d)=dt$, we obtain the claimed results.

    Finally, to prove that $c_1(\widetilde{\Hcal}(m))=t$, observe that a generator of the fiber of the Hodge line bundle over a point of $\widetilde{\Mcal}_{1,m+1}$ is given by the global differential form $d(X/Y)$, on which $\Gm$ acts with weight $1$. Here $X$ and $Y$ are inhomogeneous coordinates on the affine universal curve over $\widetilde{\Mcal}_{1,m+1}$, which is planar in the given range of $m$ (cf. \cite[Equations (1.1.3-4)]{LekiliP}).
\end{proof}

Consider the forgetful map $\oM_{1,n}(m)\dashrightarrow\oM_{1,n-1}(m)$. Suppose that the $n$-th marking lies on the core, and that the latter has level $m+1$. Then forgetting $p_n$ makes the curve $m$-unstable. The solution is to contract the core after sprouting the markings lying directly on it; this will increase the level, since there must be at least one rational tail containing at least two special points.
The problem is that the contraction is not well-defined as soon as there are at least two rational tails, because the resulting singularity has a non-trivial crimping space. In order to resolve the indeterminacy, it is enough to blow up the loci $\EB_S(m)$, where $S$ is a partition of $[n]$ into $m+1$ parts, of which at least two have $\lvert S_i\rvert\geq2$, and one is the singleton $\{n\}$; see \cite[\S 2.1]{SmythIII}.
Fortunately for us, in this paper we only need the cases $1\leq m<n\leq 4$. Even for $n=4$, there is not enough room for two rational tails and one marking on the core. So, the forgetful morphism $\oM_{1,n}(m)\dashrightarrow\oM_{1,n-1}(m)$ is well-defined in this range.

Finally, there always is a morphism $\oM_{1,m+1}(m)\to\widetilde{\Mcal}_{1,m}$. Indeed, $\oM_{1,m+1}(m)$ sits as an open inside $\widetilde{\Mcal}_{1,m+1}$, and the latter can be identified with the quotient of $\widetilde\Ucal^{sns}_{1,m+1}$ by $\Gm$, see \Cref{sec:altCompII}. Now, \cite[Proposition 1.1.5]{LekiliP} ensures that there always is a $\Gm$-equivariant morphism $\Ucal^{sns}_{1,m+1}\to \Ucal^{sns}_{1,m}$, identifying the former with the affine universal curve over the latter.

Concretely, we may think of the morphism $\oM_{1,n}(m)\to\widetilde{\Mcal}_{1,m}$ as follows. Consider the universal curve $C'$ over $\oM_{1,n}(m)$, marked with the first $m$ markings only. Rational tails and bridges made unstable by forgetting the last markings may be harmlessly contracted. 
Now, the core of $C'$ is $m$-unstable if and only if at least two of the first $m$ markings belong to a rational tail. If this is the case, we blow up the core at its markings, and then contract it (by Smyth's contraction lemma \cite[Lemma 2.13]{SmythI} or \cite[Proposition 3.7.6.1]{RSPW}). By stability, this will increase the core level, so we can reiterate until it becomes $\geq m$. For $n\geq 4$, we need to perform this operation only if $n=4$ and $m=2$. Indeed, the following two conditions must be met:
\begin{enumerate}
    \item there must be a partition $S$ having at least $m+1$ parts, and
    \item by restricting $S$ to the subset $[m]=\{1,\ldots,m\}$ we do not get the discrete partition.
\end{enumerate} 
For $n\leq 4$, this happens if and only if $n=4$ and $m=2$.
\begin{example}
For instance, consider a nodal curve with an elliptic component marked with $p_3$ and $p_4$, and a rational tail marked with $p_1$ and $p_2$; by forgetting $p_3$ and $p_4$, the elliptic core is not $2$-stable anymore. We contract it, thus producing a cusp at the attaching point of the rational tail. 
\end{example}
This is important, because contracting a genus one subcurve alters the Hodge bundle according to \Cref{prop:lambdaVSlambda}. 

Moreover observe that, whenever a core contraction takes place, the resulting singularity is an elliptic $\ell$-fold point with $\ell<m$. So, the preimage of the $m$-fold point in $\widetilde{\Mcal}_{1,m}$ remains away from the correction loci. We conclude that the normal bundle to the $m$-fold loci in $\oM_{1,n}(m)$ has the same expression as the normal bundle to the $m$-fold point in $\widetilde{\Mcal}_{1,m}$.

We summarise the previous discussion in the following:
\begin{proposition}\label{prop:normal through forgetting}
 Let $1\leq m<n\leq 4$. Forgetting the last $n-m$ markings defines a morphism $f_{n,m}\colon\oM_{1,n}(m)\to \widetilde{\Mcal}_{1,m}$. The preimage of the elliptic $m$-fold point $\Zcal$ in $\widetilde{\Mcal}_{1,m}$ is the locus of $m$-fold points in $\oM_{1,n}(m)$ such that each $p_i,\ i=1,\ldots,m,$ cleaves to a different branch of the singularity; moreover, $f_{n,m}$ is flat around these loci, and in particular their normal bundles are as in \Cref{lm:class of Z}.
\end{proposition}
\begin{corollary}\label{cor:Ell_classes}
With the notations established in the corollaries of \Cref{prop:r dagger}, we have:
 \begin{center}
\begin{tabular}{|c|c|c|}
\hline
 \diagbox{m}{n} & 3 & 4 \\ \hline
 1 & $[\mathbf{Cu}]=24\lau^2$ & $[\mathbf{Cu}]=24\lau^2$ \\
 2 & $[\mathbf{Tac}_1]+[\mathbf{Tac}_2]=24\lad^3$ & $[\mathbf{Tac}_1]+[\mathbf{Tac}_2]+[\mathbf{Tac}_{13}]+[\mathbf{Tac}_{14}]=24(\lad+\tau_{12})^3$  \\
 3 & -- & $[\mathbf{Tri}_{14}]+[\mathbf{Tri}_{24}]+[\mathbf{Tri}_{34}]=12\lat^4$ \\ \hline
\end{tabular}
\end{center}
where $\tau_{12}$ is the fundamental class of the locus of curves having a rational tail marked with $p_1$ and $p_2$.
\end{corollary}
The class $\tau_{12}$ in the corollary above appears precisely because for $n=4$ and $m=2$, in order to define the morphism to $\widetilde{\mathcal{M}}_{1,2}$, we need to contract the elliptic core of the universal curve over $\Tail{\{1,2\}}(2)$. This affects the pullback of the Hodge bundle, as explained in \Cref{prop:lambdaVSlambda}.

\section{The Chow ring of $\oM_{1,3}$}
In this section we compute the integral Chow ring of $\oM_{1,3}(m)$, for $m=2,1,0$, by a repeated application of the weighted blow-up formula \eqref{eqn:blow-up formula}.

\subsection{The Chow ring of $\oM_{1,3}(2)$ and $\oM_{1,3}(1)$}
We start from the following:
\begin{lemma}\label{prop:chow M13(2)}
    Let $\lad$ be the first Chern class of the Hodge line bundle $\HH(2)$ on $\oM_{1,3}(2)$. Then we have:
    \[A^*(\oM_{1,3}(2))=\ZZ[\lad]/(12\lad^4).\]
\end{lemma}
\begin{proof}
     By \Cref{prop:presentation M} we know that $\oM_{1,3}(2)\simeq\Pcal(1,2,2,3)$; under this isomorphism, the line bundle $\Ocal_{\Pcal}(1)$ is identified with the Hodge line bundle. The conclusion follows from the weighted projective bundle formula \cite[Theorem 3.12]{ArenaO}.
\end{proof}
We now want to apply the weighted blow-up formula \eqref{eqn:blow-up formula}, leveraging the fact that $\rho_1\colon\oM_{1,3}(1)\to\oM_{1,3}(2)$ is a weighted blow-up centered at the three points $\Tac_i$ for $i=1,2,3$ (\Cref{cor:picture M13}). As explained in \S \ref{sec:chow wb}, for this we need to compute the following:
\begin{enumerate}
    \item the fundamental class of the center, i.e. the three points $\Tac_{i}$;
    \item the $\Gm$-equivariant top Chern class of the normal bundle of $\Tac_{i}$;
    \item surjectivity and the kernel of the pullback homomorphism $A^*(\oM_{1,3}(2))\to A^*(\Tac_{i})$.
\end{enumerate}
By Corollary \ref{cor:Ell_classes} we know that $[\Tac_1]+[\Tac_2]=24\lad^3$, hence by symmetry each $\Tac_i$ has class $12\lad^3$.
Point (3) is straightforward since $\Tac_i$ is a point for every $i$ (\Cref{cor:picture M13}). By the same token, the normal bundle is trivial with $\Gm$-action of weights $(2,3,4)$:
\begin{equation}\label{eq:top Chern Tac}
    c_3^{\Gm}(\mathcal{N}_{\Tac_i}) = 24t^3.
\end{equation}

Let us proceed by blowing up one tacnodal point at a time: We start by computing the Chow ring of the weighted blow-up $\oM_{1,3}(2)'$ of $\oM_{1,3}(2)$ at the point $\Tac_1$. The pullback homomorphism $A^*(\oM_{1,3}(2))\to A^*(\Tac_1)$ is surjective, and its kernel is generated by $\lad$. We may hence apply the simplified weighted blow-up formula \eqref{eqn:blow-up formula}:
\begin{align*}
    A^*(\oM_{1,3}(2)') &= A^*(\oM_{1,3}(2))[\delta_1]/(\lad\delta_1,12\lad^3-24\delta_1^3).
\end{align*}
where we denoted by $\delta_1$ the fundamental class of the exceptional divisor, representing curves with a node separating $p_1$ on the core from $p_2$ and $p_3$ on a rational tail.
\begin{remark}
    Notice the sign in the last relation, which is obtained by evaluating the equivariant top Chern class of the normal bundle of the blow-up centre by setting the equivariant parameter $t$ of \eqref{eq:top Chern Tac} equal to the negative of the fundamental class of the exceptional divisor.
\end{remark}
When we blow up $\Tac_2$ (respectively $\Tac_3$) we get a similar result, with the only difference that now the kernel of the pullback homomorphism is generated by $\lad$ and $\delta_1$ (respectively $\lad$, $\delta_1$ and $\delta_2$). In the end we obtain the following:
\begin{lemma}\label{prop:chow M131 1} The Chow ring of $\oM_{1,3}(1)$ can be written as follows:
    \[ A^*(\oM_{1,3}(1))=\ZZ[\lad,\de_1,\de_2,\de_3]/(12\lad^4,\{\de_i\de_j\}_{i\neq j},\{\lad\de_i,-24\de_i^3+12\lad^3\}_{i=1,2,3}). \]
\end{lemma}
\begin{remark}
 The first relation is now redundant, as for any $i$ we can write
 \[ 12\lad^4 = \lad (-24\delta_i^3+12\lad^3) + 24\delta_i^2 (\delta_i\lad).\]
\end{remark}

On $\oM_{1,3}(1)$ there are two families of curves, the universal $1$-stable curve, and the $2$-stable curve pulled back from $\oM_{1,3}(2)$. Their Hodge line bundles are related by \Cref{prop:lambdaVSlambda}. Accordingly, we make a change of variable: $$\lau=\lad-\delta_1-\delta_2-\delta_3,$$ where $\lau$ is the first Chern class of the Hodge line bundle $\HH(1)$ of the $1$-stable curve. Moreover, recall that $\oM_{1,3}(1)$ is the same as Schubert's stack of pseudo-stable curves.
\begin{proposition}\label{prop:chow M13(1)}
The Chow ring of the moduli stack of pseudo-stable $3$-pointed curves of genus one is:
 \begin{equation*}
 A^*(\oM_{1,3}^{ps})\simeq\ZZ[\lau,\de_1,\de_2,\de_3]/I^{ps}_{1,3}
\end{equation*}
where 
the ideal of relations $I^{ps}_{1,3}$ is generated by:
\begin{align*}
    &\de_i^2+\lau\de_i=0 &\text{for }i=1,2,3,\\
    &\de_i\de_j=0 &\text{for }i\neq j,\\
    &24(\de_i^3-\de_j^3)=0 &\text{for }i\neq j,\\
    &12(\lau^3+\de_1^3+\de_2^3-\de_3^3)=0.
\end{align*}
\end{proposition}

\subsection{The Chow ring of $\oM_{1,3}$}
Recall from \Cref{cor:picture M13} that $\rho\colon\oM_{1,3}\to\oM_{1,3}(1)$ is a weighted blow-up of weights $(4,6)$ centred in the locus of cuspidal curves $\Cu$. The exceptional divisor can be identified with the locus of unmarked elliptic tails $\Delta_{\emptyset}$.

In order to apply the simplified weighted blow-up formula, we need the following ingredients:
\begin{enumerate}
    \item the surjectivity and kernel of the pullback homomorphism $A^*(\oM_{1,3}(1)) \to A^*(\Cu)$;
    \item the fundamental class of $\Cu$;
    \item the $\Gm$-equivariant top Chern class of the normal bundle of $\Cu$.
\end{enumerate}
\subsubsection{Pinching and pulling}
We aim to study the Gysin pullback:
\[ \iota^*\colon A^*(\oM_{1,3}(1))\to A^*(\mathbf{Cu}). \]
Pinching (collapsing the tangent line) at the marking $\star$ identifies $\Cu$ with $\oM_{0,3\sqcup\star}$, which in turn is isomorphic to $\PP^1$. Hence, 
\[A^*(\mathbf{Cu})=A^*(\oM_{0,3\sqcup\star})=A^*(\PP^1)=\ZZ[h]/(h^2),\]
where any boundary point and any $\psi$-class on $\oM_{0,3\sqcup\star}$ is identified with the generator $h$ (see \cite[\S 1.5]{Kock}). 

The boundary divisor $\delta_i$ intersects $\Cu$ transversely in one point, representing a cusp marked by $p_i$ with a rational tail marked by $p_j$ and $p_k$: its class can therefore be identified with the generator of the Chow ring $h$; in particular, the pullback $\iota^*$ is surjective. On the other hand, Proposition \ref{prop:SmythLambda} gives us $\iota^*\lau= -\psi_{\star} = -h$.
Putting everything together, we obtain the following, which is ingredient (1) for applying the weighted blow-up formula.
\begin{lemma}\label{lemma:iota surj}
    The pullback homomorphism $\iota^*:A^*(\oM_{1,3}(1))\to A^*(\Cu)$ is surjective and
    \[ \ker(\iota^*) = (\lau+\de_1,\lau+\de_2,\lau+\de_3). \]
\end{lemma}

\subsubsection{End of the computation}
The normal bundle of $\Cu$ can be pulled back along $f\colon \oM_{1,3}(1) \to \widetilde{\mathcal{M}}_{1,1}$, as explained in \Cref{prop:normal through forgetting}, and by \Cref{lm:class of Z} it is equal to $$f^*( \widetilde{\mathcal{H}}(1)^{\otimes 4}\oplus \widetilde{\mathcal{H}}(1)^{\otimes 6})|_{\Cu}=\mathcal{H}(1)^{\otimes 4} \oplus \mathcal{H}(1)^{\otimes 6}|_{\Cu}.$$ This provides ingredient (3). Finally, ingredient (2) is provided to us by Corollary \ref{cor:Ell_classes}. After a straightforward computation, we obtain
\begin{equation*}
 A^*(\oM_{1,3})=\ZZ[\lambda_1,\de_1,\de_2,\de_3,\de_{\emptyset}]/ I'_{1,3}
\end{equation*}
where $I'_{1,3}$ is generated by the relations
\begin{align*}
    &\de_i\de_j&\text{ for }i,j\in \{1,2,3\},\\
    &\de_i^2+\lau\de_i&\text{ for }i\in\{1,2,3\},\\
    &12\lau^2(\lau+\de_i+\de_j-\de_k)&\text{ for }i,j,k\in\{1,2,3\},\\
    &\de_{\emptyset}(\lau+\de_i)&\text{ for }i\in\{1,2,3\},\\
    &24(\lau-\de_{\emptyset})^2.&
\end{align*}
\Cref{prop:lambdaVSlambda} gives us the change of variables $\la=\lau-\de_{\emptyset}$, under which we get:
\begin{theorem}\label{thm:chow M13}
    The integral Chow ring of the moduli stack of $3$-marked stable curves of genus $1$ is
    \begin{equation*}
    A^*(\oM_{1,3})=\ZZ[\la,\de_1,\de_2,\de_3,\de_{\emptyset}]/ I_{1,3}
    \end{equation*}
   where 
   $I_{1,3}$ is generated by ten relations in degree $2$:
   \begin{subequations}
    \begin{align}
    &24\la^2=0 & \subtag{1}\\
    &\de_{\emptyset}^2=-\la\de_{\emptyset}-\de_3\de_{\emptyset} &\subtag{2}\\
    &\de_i^2=-\la\de_i-\de_3\de_{\emptyset}, & i=1,2,3 \subtag{2+i}\\
    &\de_i\de_j=0 \subtag{5+k}, & \{i,j,k\}=\{1,2,3\}\\
    &\de_{\emptyset}\de_1=\de_{\emptyset}\de_j, & j=2,3 \subtag{7+j},
    \end{align}
    and by one extra (torsion) relation in degree $3$, namely
    \begin{align}
    &12(\la+\de_{\emptyset})^2(\la+\de_{\emptyset}+\de_i+\de_j-\de_k)=0. & \subtag{11}
    \end{align}
\end{subequations}
\end{theorem}

\begin{remark}
 Tensoring with $\mathbb{Q}$, we recover \cite[Theorem 3.3.2]{Belorousski}.
\end{remark}

\section{The Chow ring of $\oM_{1,4}$}
We compute the integral Chow ring of $\oM_{1,4}$ by the same strategy. The main novelty is that the birational map $\oM_{1,4}(1)\dashrightarrow\oM_{1,4}(2)$ is not regular, hence we have to climb down a blow-up.
\subsection{The Chow ring of $\oM_{1,4}(3)$ and $\oM_{1,4}(2)$}
Our starting point is the Chow ring of $\oM_{1,4}(3)$, which can be easily computed using the fact that $\oM_{1,4}(3)\simeq\Pcal(1,1,1,2,2)$ (\Cref{prop:presentation M} and the formula for the Chow ring of a weighted projective stack \cite[Theorem 3.12]{ArenaO}).
\begin{lemma}\label{prop:Chow M14(3)}
    We have
     \[A^*(\oM_{1,4}(3))=\ZZ[\lat]/(4\lat^5)\]
     where $\lat$ is the first Chern class of the Hodge line bundle $\Hcal(3)$.
\end{lemma}
To compute the Chow ring of $\oM_{1.4}(2)$, we leverage \Cref{cor:picture M14}: there is a morphism
\[\rho_2\colon \oM_{1,4}(2) \longrightarrow \oM_{1,4}(3) \]
which is a weighted blow-up of weights $2$, $3$ and $4$ centered at the locus of elliptic curves with a triple point. To apply the weighted blow-up formula we need as usual the following ingredients:
\begin{enumerate}
     \item the surjectivity and kernel of the pullback homomorphism $A^*(\oM_{1,4}(3)) \to A^*(\Tri_{ij})$;
    \item the fundamental class of $\Tri_{ij}$;
    \item the $\Gm$-equivariant top Chern class of the normal bundle of $\Tri_{ij}$.
\end{enumerate}
The second ingredient is given to us by Corollary \ref{cor:Ell_classes}: the flat forgetful morphism $\oM_{1,4}(3)\to \widetilde{\Mcal}_{1,3}$ pulls $[\Zcal]=12t^4$ back to $[\Tri_{14}]+[\Tri_{24}]+[\Tri_{34}]=12\lat^4$, so by the action of the symmetric group
\[ [\Tri_{ij}] = 4\lat^4 \text{ for all }\{i,j\}\subseteq[4]. \]
Being points, the normal bundle is trivial (with $\Gm$-action of weights $(2,3,4)$), and restriction is surjective on Chow groups. From the blow-up formula \eqref{eqn:blow-up formula} 
we get:
\begin{equation}\label{eq:chow M14(2) pre}
A^*(\oM_{1,4}(2))=\ZZ[\lat,\{\tau_{ij}\}_{\{i,j\}\subset [4]}]/I_{1,4}(2)',
\end{equation}
where $\tau_{ij}$ is the class of the exceptional divisor $\TT_{ij}$, which can be identified with the locus of curves with a node separating the markings $h,k$ on the core from the markings $i,j$ on a rational tail\footnote{We reserve our usual notation $\Delta_{hk}$ for the analogous locus in $\oM_{1,4}$ and $\oM_{1,4}(1)$; as we shall see in \Cref{prop:gen M141}, the classes $\tau_{ij}$ and $\delta_{hk}$ do not quite agree.}, 
and $I_{1,4}(2)'$ is generated by
\begin{subequations}
    \begin{align}
    &\lat\tau_{ij}=0, & \{i,j\}\subset [4], \subtag{1--6}\\
    &\tau_{ij}\tau_{hk}=0, & \{i,j \}\neq \{h,k\},  \subtag{7--21}\\
    &12\tau_{ij}^4+4\lat^4=0. & \subtag{22--27}
    \end{align}
\end{subequations}
The relation $4\lat^5$ is redundant and we thus omit it. 

\Cref{prop:lambdaVSlambda} yields the change of variables $\lat=\lad+\sum_{1\leq i<j\leq 4}\tau_{ij}$.
After performing this substitution, we obtain the following.
\begin{proposition}\label{prop:chow M14(2)}
The Chow ring of the moduli stack of $4$-marked $2$-stable curves of genus one is:
 \begin{equation*}
 A^*(\oM_{1,4}(2))=\ZZ\left[\lad,\{\tau_{ij}\}_{\{i,j\}\subset [4]}\right]/I_{1,4}(2)
\end{equation*}
where the ideal of relations $I_{1,4}(2)$ is generated by:
\begin{subequations}
\begin{align}
    \lad\tau_{ij}+\tau_{ij}^2, &&\{i,j\}\subset [4], \subtag{1--6}\\
    \tau_{ij}\tau_{hk}, &&\{ i,j \}\neq \{h,k\}, \subtag{7--21}\\
    12\tau_{ij}^4 + 4(\lad^4-\sum_{r,s}\tau_{rs}^4) && \subtag{22--27}
\end{align}
\end{subequations}
The rational Chow ring of the moduli space $\overline{M}_{1,4}(2)$ has Hilbert series $1+7t+7t^2+7t^3+t^4$.
\end{proposition}

\subsection{The Chow ring of $\oM_{1,4}(1)$}
We distinguish the tacnodal loci in $\oM_{1,4}(2)$ into two types: those separating one marking $i$ from the other three (there are four of them, call them $\Tac_i$), and those separating two markings $i,j$ from the other two (there are three, call them $\Tac_{ij}$). The birational map $\oM_{1,4}(1)\dashrightarrow \oM_{1,4}(2)$ can be resolved by a suitable blow-up (\Cref{cor:picture M14}):
    \bcd
& \oM_{1,4}(\frac{3}{2})\ar[dl,"\rho_{\frac{3}{2}}"]\ar[dr,"\rho_2"]& \\
\oM_{1,4}(1)\ar[rr,dashed] & & \oM_{1,4}(2)
    \ecd
Here $\rho_2$ is a blow-up of weights $(2,3,4)$ in the tacnodal locus. We denote the classes of the exceptional divisors  by $\beta$ in \Cref{prop:3/2} below. The morphism $\rho_{\frac{3}{2}}$ is an ordinary blow-up of the elliptic bridges $\EB_{\{ij\},\{hk\}}$, whose exceptional divisors can be identified with those over $\mathbf{Tac}_{ij}$; the loci $\EB_{\{ij\},\{hk\}}$ are thus flipped. In particular, $\beta_{ij}$ pushes forward to $0$ along $\rho_{\frac{3}{2}}$, while $\beta_i$ does not.

We will first compute the Chow ring of $\oM_{1,4}(\frac{3}{2})$ via the weighted blow-up formula, and then identify the Chow ring of $\oM_{1,4}(1)$ as a subring of the former.
\subsubsection{The Chow ring of $\oM_{1,4}(\frac{3}{2})$}
We need the following ingredients:
\begin{enumerate}
    \item the surjectivity and kernel of the pullback homomorphism to $A^*(\Tac_{ij})$ and $A^*(\Tac_{i})$;
    \item the fundamental classes of $\Tac_{ij}$ and $\Tac_i$;
    \item the $\Gm$-equivariant top Chern class of the normal bundle of $\Tac_{ij}$ and $\Tac_i$.
\end{enumerate}

For the first point, we need to compute the pullback of the generators of $A^*(\oM_{1,4}(2))$ in $A^*(\Tac_i)$ (respectively $A^*(\Tac_{ij})$).
Let $\varphi_i\colon \PP^1\simeq\oM_{0,\{\star,h,j,k\}}\to \oM_{1,4}(2)$ be the pinching morphism, defined as follows: first, consider the auxiliary stack $\mathcal{X}$ whose objects are triples $(C\to S,D\to S,\alpha)$, where
\begin{enumerate}
    \item[(i)] $C\to S$ is a stable curve of genus zero with markings indexed by $\star$, $h$, $j$ and $k$,
    \item[(ii)] $D\to S$ is a smooth genus zero curve with two markings indexed by $\star'$ and $i$, and
    \item[(iii)] $\alpha\colon (\mathfrak{m}_{\star}/\mathfrak{m}^2_{\star})^\vee \simeq (\mathfrak{m}_{\star'}/\mathfrak{m}^2_{\star'})^\vee$ is an isomorphism.
\end{enumerate}
There is a natural morphism $\mathcal{X} \to \oM_{1,4}(2)$ whose image is $\Tac_i$: this morphism is given by attaching $C\to S$ and $D\to S$ at the markings indexed by $\star$ and $\star'$ by gluing the two tangent spaces via $\alpha$.
On the other hand, the projection $\mathcal{X} \to \oM_{0,\{\star, h,j,k\}}$ is an isomorphism: indeed, the fibres consist of the choice of a smooth $2$-marked genus zero curve $D\to S$ and an isomorphism $\alpha$. 
There is a unique isomorphism class of such $D$s, with automorphism group isomorphic to $\Gm$. The latter acts with weight one on the tangent space at $\star'$, and hence it acts freely and transitively on the set of isomorphisms $\alpha$. 
By composing the inverse of the projection with the gluing map, we obtain the morphism $\varphi_i\colon\oM_{0,\{\star,h,j,k\}}\to \oM_{1,4}(2)$ with image $\Tac_i$.

Similarly, let $\varphi_{ij}\colon\PP^1\to \oM_{1,4}(2)$  with image $\Tac_{ij}$ be the crimping morphism: the pointed normalisations of these tacnodal curves have no moduli, since they consist of two $3$-pointed rational curves $C\to S$ and $D\to S$, marked with $\{i,j,\star\}$ and $\{h,k,\star'\}$ respectively. Let $\PP^1=\PP(T_{\star}C \oplus T_{\star'}D)$ be the compactified moduli space of attaching data: the points away from $0$ and $\infty$ correspond to the choice of an isomorphism $(\mathfrak{m}_{{\star}}/\mathfrak{m}_{{\star}}^2)^\vee \simeq (\mathfrak{m}_{{\star'}}/\mathfrak{m}_{{\star'}}^2)^\vee$;
over $0$ and $\infty$, a rational tail is sprouted, supporting the markings $\{i,j\}$ and $\{h,k\}$ respectively, so the automorphism group of the curve acts transitively on the attaching data, as in the previous paragraph. More precisely, we blow up $C\times\PP^1$ at $(\star,0)$ and $D\times\PP^1$ at $(\star',\infty)$, and then glue the proper transforms of $\star$ and $\star'$ according to their common projection in $\PP^1$. 

Therefore, since both $\Tac_i$ and $\Tac_{ij}$ are isomorphic to $\PP^1$, their Chow ring can be identified with $\ZZ[h]/(h^2)$.
\begin{lemma}\label{lm:other coeff Tac}
    Fix $1\leq i<j\leq 4$, and let $\TT_{ij}$ denote the locus of curves with a rational tail marked by $i,j$ inside $\oM_{1,4}(2)$. 
    Then:
    \begin{enumerate}
        \item $\int_{\Tac_\bullet}\varphi_\bullet^*\lad=-1$ for both $\bullet=i$ and $\bullet=ij$;
        \item for $1\leq h<k\leq 4$, $\Tac_{ij}$ intersects the divisor $\TT_{hk}$ if and only if $\{h,k\}=\{i,j\}$ or $[4]\setminus\{i,j\}$, in which cases they intersect transversely at one point; 
        \item $\Tac_i$ intersects the divisor $\TT_{hk}$ if and only if $\{h,k\}\subset[4]\setminus\{i\}$, in which cases they intersect transversely at one point. 
    \end{enumerate}
    In particular, the kernel of the Gysin pullback to $\Tac_i$ and $\Tac_{ij}$ is determined as follows:
    \begin{align*}
    \ker(A^*(\oM_{1,4}(2)) \longrightarrow A^*(\Tac_{i})) &= (\{\tau_{ij}\}_{j\neq i}, \{\lad + \tau_{hk}\}_{h,k\neq i},\lad^2 ). \\
    \ker(A^*(\oM_{1,4}(2)) \longrightarrow A^*(\Tac_{ij})) &= (\{\tau_{ih}\}_{h\neq j}, \{\tau_{jk}\}_{k\neq i},\lad + \tau_{hk}, \lad + \tau_{ij},\lad^2 ).
\end{align*}
\end{lemma}
\begin{proof}
    Point (1) is Proposition \ref{prop:SmythLambda}, the others follow from the deformation theory of nodes.
\end{proof}

Now, we want to determine the classes of $\Tac_i$ and $\Tac_{ij}$ in $A^3(\oM_{1,4}(2))$ in terms of $\lat$ and $\tau$'s. 
From Corollary \ref{cor:Ell_classes} we know that:
\begin{equation}\label{eqn:TacLoci}
    [\Tac_1] + [\Tac_2] +  [\Tac_{13}] + [\Tac_{14}] = f_{4,2}^*[\Zcal] = 24(\lad+\tau_{12})^3=24\lad^3+24\tau_{12}^3.
\end{equation}
We exploit the $S_4$-action (by permuting the markings) to simplify the intersection-theoretic calculations. As an $S_4=S(1,2,3,4)$-module, $A^3(\oM_{1,4}(2))$ can be identified with $\on{triv}(\ell)\oplus\Lambda^2\mathbb C^4_{\text{std}}$, the former with basis $\lad^3=:\ell$, and the latter with basis $\tau_{ij}^3=:e_{ij}$. Note that $[\Tac_1]$ is invariant under $S_3=S(2,3,4)$, hence its class can be written as:
\[[\Tac_1]=u_1\ell+a_1(e_{12}+e_{13}+e_{14})+b_1(e_{23}+e_{24}+e_{34}).\]
Similarly
\[[\Tac_2]=u_2\ell+a_2(e_{12}+e_{23}+e_{24})+b_2(e_{13}+e_{14}+e_{14}).\]
Note that $[\Tac_{13}]$ is invariant under $S_2\times S_2=S(1,3)\times S(2,4)$, hence its class can be written as:
\[[\Tac_{13}]=v_1\ell+c_1e_{13}+d_1e_{24}+f_1(e_{12}+e_{14}+e_{23}+e_{34}).\]
Similarly
\[[\Tac_{14}]=v_2\ell+c_2e_{14}+d_2e_{23}+f_2(e_{12}+e_{13}+e_{24}+e_{34}).\]
By the $S_4$-action we see that $$u_1=u_2=:u,v_1=v_2=:v,a_1=a_2=:a, b_1=b_2=:b,c_1=d_1=c_2=d_2=:c,f_1=f_2=:f.$$ Grouping the coefficients of $e_{ij}$, Equation \eqref{eqn:TacLoci} reduces to the following linear system:
\[\begin{pmatrix}
2 & 0 & 0 & 2 \\
1 & 1 & 1 & 1 \\
0 & 2 & 0 & 2
\end{pmatrix}
\begin{pmatrix} a\\b\\c\\f\end{pmatrix}=
\begin{pmatrix} 24\\0\\0\end{pmatrix}\]
Furthermore, by intersecting $\Tac_i$ with $\tau_{ij}$ we obtain $u=a$: indeed, we know from \Cref{lm:other coeff Tac} that $\tau_{ij}$ is in the kernel of the pullback to $\Tac_i$, hence $[\Tac_i]\cdot \tau_{ij}=0$; on the other hand, by \Cref{prop:chow M14(2)} we have $\tau_{ij}\tau_{hk}=0$ unless $\{i,j\}=\{h,k\}$, and $\tau_{ij}\lad=-\tau_{ij}^2$. Putting everything together, we deduce that $[\Tac_i]\cdot \tau_{ij}=(a-u)\tau_{ij}^4$, hence $a-u=0$ because $A^4(\oM_{1,4}(2))\otimes \mathbb{Q}\simeq \mathbb{Q}\cdot \tau_{ij}^4$.  

With a similar argument, by intersecting $\Tac_{ij}$  with $\tau_{ik}$ we obtain $v=f$. We need one more equation, which we are going to obtain by computing the pullback of $\lat$ from $\oM_{1,4}(3)$ to $\Tac_i$: this will determine $u$ quite directly. Recall the morphism $\rho_3\colon\oM_{1,4}(2)\to\oM_{1,4}(3)$.
\begin{lemma}\label{lem:tacnodalloci1}
     The composition $\rho_3\circ \varphi_i:\Tac_i\simeq\oM_{0,4}\to\oM_{1,4}(3)$ maps the three special points to the $3$-elliptic loci $\Tri_{jk}$, $i\notin\{j,k\}\subset[4]$. The following hold:
    \begin{align} \label{pippo}
        \int_{\Tac_i}\varphi_i^*\psi_j&=1,\\\label{pluto}
        \int_{\Tac_i}\varphi_i^*\lad&=-1,\\\label{topolino}        \int_{\Tac_i}\varphi_i^*\rho_3^*\lat&=2.
    \end{align}
    Therefore, $u=8$ and $v=4$.
\end{lemma}
\begin{proof}
    \Cref{pippo} is a genus zero calculation, since crimping happens away from the marking.
    
    \Cref{pluto} follows from Proposition \ref{prop:SmythLambda} and a genus zero computation.\footnote{Alternatively, from Proposition \ref{prop:psiVSlambda} we obtain $\int_{\Tac_i}\varphi_i^*\lad=\int_{\Tac_i}\varphi_i^*\psi_j-2$, since the universal $2$-stable curve over $\Tac_i$ is obtained by gluing a trivial family $(\PP^1,i,\star)$ with a moving $(\PP^1,h,j,k,\star)$, so the marking $j$ ends up on a rational tail twice. We conclude by \Cref{pippo}.}

    From Proposition \ref{prop:lambdaVSlambda} it follows that $\int_{\Tac_i}\varphi_i^*\rho_3^*\lat=\int_{\Tac_i}\varphi_i^*\lad+3$, since the universal $3$-stable family is obtained by sprouting the lonely markings and contracting the tacnode over the three special points. In view of this, \Cref{topolino} follows from \Cref{pluto}.

    We obtain $2u+2v=24$ by looking at the coefficients of $\lad^3$ in \Cref{eqn:TacLoci}, hence $v$ is determined by $u$. Now observe that, because of the relations $\lambda_{(3)}\tau_{ij}=0$ (see \eqref{eq:chow M14(2) pre}) and $\lat=\lad+\sum \tau_{ij}$ (see \Cref{prop:lambdaVSlambda}), the projection formula gives
    \[\rho_{3,*}(\rho_3^*\lat \cdot [\Tac_{i}]) = u\lat^4.\]
    On the other hand, we deduce from \Cref{topolino} that
    \[\rho_{3,*}(\rho_3^*\lat \cdot [\Tac_{i}]) = \rho_{3,*}\rho_3^*\varphi_{i,*}\varphi_i^*\lat = \rho_{3,*}\varphi_{i,*}(2[\text{pt}]) = 2[\Tri_{jk}] = 8\lat^4.\]
\end{proof}
We may now solve the linear system above and find $u=a=-c=8,\ v=-b=f=4$.
\begin{proposition}
    The classes of the tacnodal loci are
    \begin{align*}
    [\Tac_i]&=8\lad^3-4(\tau_{hj}^3+\tau_{hk}^3+\tau_{jk}^3) +8(\tau_{ih}^3+\tau_{ij}^3+\tau_{ik}^3),\\
    [\Tac_{ij}]&=4\lad^3-8(\tau_{ij}^3+\tau_{hk}^3)+4(\tau_{ih}^3+\tau_{ik}^3+\tau_{jh}^3+\tau_{jk}^3).
\end{align*}
\end{proposition}
The equivariant normal bundle is given by \Cref{lm:class of Z}.
We can now compute the Chow ring of $\oM_{1,4}(\frac{3}{2})$ by applying the weighted blow-up formula \eqref{eqn:blow-up formula}.

\begin{proposition}
The Chow ring of the moduli stack of \emph{aligned} $4$-marked $1$-stable curves of genus $1$ is
generated by the symbols $\lad,\tau_{ij},\beta_i,\beta_{ij}=\beta_{hk}$ (where $[4]=\{i,j,h,k\}$), subject to the following relations:
\begin{align*}
    &(1-6)  &&\lad\tau_{ij}+\tau_{ij}^2, &&\{i,j\}\subset [4],\\
    &(7-21) &&\tau_{ij}\tau_{hk}, &&\{h,k\}\neq \{i,j\},\\
    &(22-27) &&12\tau_{ij}^4 + 4(\lad^4-\sum_{r,l}\tau_{rl}^4),&& \\
    &(28-39) &&\beta_i\tau_{ij} &&j\in[4]\setminus\{i\},\\
    &(40-51) &&\beta_i(\lad+\tau_{jk}), &&j,k\in[4]\setminus\{i\}, \\
    &(52-57) &&\beta_i\beta_j, &&j\in[4]\setminus\{i\}, \\
    &(58-61) &&24\beta_i(3\lad\beta_i-\beta_i^2)+[\Tac_i],&& \\
    & (62-67) &&\beta_{ij}\tau_{ih}, &&h\in[4]\setminus\{i,j\}, \\
    & (68-73) &&\beta_{ij}(\lad+\tau_{ij}), && \\
    & (74-85) &&\beta_{ij}\beta_h, &&h\in[4], \\
    & (86-91) &&\beta_{ij}\beta_{hk},&&\{h,k\}\neq\{i,j\},\\
    & (92-94) &&24\beta_{ij}(3\lad\beta_{ij}-\beta_{ij}^2)+[\Tac_{ij}]. &&
\end{align*}
The Hilbert series of the rational Chow ring of the moduli space $\overline{M}_{1,4}(\frac{3}{2})$ is equal to $1+14t+21t^2+14t^3+t^4$.
\end{proposition}

According to Proposition \ref{prop:lambdaVSlambda}, we make the substitution $\lau=\lad-\sum_i\beta_i-\sum_{i,j}\beta_{ij}$.

\begin{proposition}\label{prop:3/2}
The Chow ring of $\oM_{1,4}(\frac{3}{2})$ is
generated by $\lau,\tau_{ij},\beta_i,\beta_{ij}=\beta_{hk}$ with:
\begin{align*}
    &(1-6)  &&\tau_{ij}(\lau+\beta_h+\beta_k+\beta_{ij}+\tau_{ij}), &&\{i,j\}\subset [4],\\
    &(7-21) &&\tau_{ij}\tau_{hk}, &&\{h,k\}\neq \{i,j\},\\
    &(22-27) &&12\tau_{ij}^4 + 4((\lau+\sum_i\beta_i+\sum_{i,j}\beta_{ij})^4-\sum_{r,l}\tau_{rl}^4),&& \\
    &(28-39) &&\beta_i\tau_{ij} &&j\in[4]\setminus\{i\},\\
    &(40-51) &&\beta_i(\lau+\beta_i+\tau_{jk}), &&j,k\in[4]\setminus\{i\}, \\
    &(52-57) &&\beta_i\beta_j, &&j\in[4]\setminus\{i\}, \\
    &(58-61) &&24\beta_i(3\lau\beta_i+2\beta_i^2)+[\Tac_i],&& \\
    & (62-67) &&\beta_{ij}\tau_{ih}, &&h\in[4]\setminus\{i,j\}, \\
    & (68-73) &&\beta_{ij}(\lau+\beta_{ij}+\tau_{ij}), && \\
    & (74-85) &&\beta_{ij}\beta_h, &&h\in[4], \\
    & (86-91) &&\beta_{ij}\beta_{hk},&&\{h,k\}\neq\{i,j\},\\
    & (92-94) &&24\beta_{ij}(3\lau\beta_{ij}+2\beta_{ij}^2)+[\Tac_{ij}]. &&
\end{align*}
\end{proposition}

\subsubsection{The Chow ring of $\oM_{1,4}(1)$}

The pullback homomorphism
\[ \rho_{\frac{3}{2}}^*:A^*\left(\oM_{1,4}(1)\right) \longrightarrow A^*\left(\oM_{1,4}\left(\frac{3}{2}\right)\right) \]
is injective \cite[Proposition 6.7(b)]{Fulton}, hence if we know a set of generators of $A^*(\oM_{1,4}(1))$, we can deduce a full presentation by looking at the subring of $A^*(\oM_{1,4}(\frac{3}{2}))$ generated by their pullback.
We recall the following useful lemma.
\begin{lemma}\label{pushpull}
 Let $\iota\colon Z\hookrightarrow X$ be a regularly embedded subvariety of a smooth ambient space. If Gysin pullback along $\iota$ is surjective, then $\on{im}(\iota_*)$ is generated as an ideal by the fundamental class $[Z]$.
\end{lemma}
\begin{proof}
 If $\alpha=\iota^*\widetilde\alpha$, the projection formula implies that $\iota_*\alpha=\widetilde\alpha\cdot[Z]$.
\end{proof}

We recall our running convention: for $S\subseteq[4]$, the boundary divisor $\Delta(1)_S\subset \oM_{1,4}(1)$ is the substack of curves with a separating node such that the markings $p_i$, $i\in S$, lie on the genus one side. Its fundamental class is denoted by $\delta_{S}$. On the other hand, we have denoted by $\TT_{ij}$ (and $\tau_{ij}$) the locus (and its class) of rational tails marked by $i$ and $j$ in $\oM_{1,4}(2)$.
\begin{proposition}\label{prop:gen M141}
    The Chow ring of $\oM_{1,4}(1)$ is generated 
    by $\lambda_{(1)}$, $\{\delta_{i}\}_{i\in [4]}$ and $\{\delta_{ij}\}_{\{i,j\}\subset [4]}$. We have:
    $$\rho_{\frac{3}{2}}^*(\lau) = \lau,\quad \rho_{\frac{3}{2}}^*(\delta_i)=\beta_i, \quad \rho_{\frac{3}{2}}^*(\delta_{ij})=\tau_{hk}+\beta_{hk}.$$
\end{proposition}
\begin{proof}
    The second statement is a simple computation of the total transform of $\Delta(1)_i$ and $\Delta(1)_{ij}$.
    
    For proving the first statement, let $\BB(\frac{3}{2})=\cup \BB_{ij}$ be the exceptional divisor of $\rho_{\frac{3}{2}}$, where $\BB_{ij}$ is the exceptional divisor of $\rho_{2}$ that lies over $\Tac_{ij}$. We claim that $A^*(\oM_{1,4}(\frac{3}{2})\smallsetminus \BB(\frac{3}{2}))$ is generated by $\lau$, $\{\beta_i\}_{i\in [4]}$ and $\{\tau_{ij}\}_{\{i,j\}\subset [4]}$. This claim follows from the previous lemma and the surjectivity of the pullback homomorphism for every $i,j$:
    $$\iota^*_{ij}:A^*\left(\oM_{1,4}\left(\frac{3}{2}\right)\right)\longrightarrow A^*(\BB_{ij}).$$
   Recall the isomorphism:
    $$\BB_{ij}\simeq \oM_{1,2}(1)\times\PP^1\simeq \Pcal(2,3,4)\times\PP^1,$$ 
    So, the Chow ring of $\BB_{ij}$ is generated (as a ring) by the two hyperplane classes $h_{\Pcal}$ and $h_{\PP^1}$. We have $\iota_{ij}^*\beta_{ij}=-h_{\PP^1}$ and $\iota_{ij}^*(\rho_{\frac{3}{2}}^*\lau) = h_{\Pcal}$, so $\iota_{ij}^*$ is surjective.

    Observe that $\rho_{\frac{3}{2}}$ induces an isomorphism of $\oM_{1,4}(\frac{3}{2})\smallsetminus \BB(\frac{3}{2})$ with $\oM_{1,4}(1)\smallsetminus\mathbf{EB}(1)$, where $\mathbf{EB}(1)=\cup\mathbf{EB}_{ij}$ is the codimension two substack of elliptic bridges, and $\mathbf{EB}_{ij}\simeq \Pcal(2,3,4)$, whose Chow ring is generated by the restriction of $\lau$. Moreover, the locus $\mathbf{EB}_{ij}$ is the transversal intersection of the divisors $\Delta(1)_{ij}$ and $\Delta(1)_{hk}$, thus $[\mathbf{EB}_{ij}]=\delta_{ij}\delta_{hk}$. We conclude that 
    $A^*(\oM_{1,4}(1))$ is generated by $\lau$, $\{\delta_{i}\}_{i\in [4]}$, $\{\delta_{ij}\}_{\{i,j\}\subset [4]}$ as a ring.
\end{proof}
\begin{corollary}
   $A^*(\oM_{1,4}(1))$ is generated by $\lau$, $\{\beta_{i}\}_{i\in [4]}$, $\{\beta_{ij}+\tau_{ij}\}_{\{i,j\}\subset [4]}$ as a subring of $A^*(\oM_{1,4}(\frac{3}{2}))$.
\end{corollary}
We find a presentation of $A^*(\oM_{1,4}(1))$ by eliminating the remaining variables from $A^*(\oM_{1,4}(\frac{3}{2}))$. The following has been produced by Sage \cite{sagemath}; we do not know if it is a minimal presentation.

\begin{proposition}\label{prop:chow M14(1)}
The Chow ring of the moduli stack of $4$-marked pseudo-stable curves of genus one is:
\[A^*(\oM_{1,4}^{ps})=\ZZ[\lau,\{\delta_{ij}\}_{\{i,j\}\subset [4]},\{\delta_i\}_{i\in [4]}]/I_{1,4}^{ps}\]
where $I_{1,4}^{ps}$ is generated by the following relations:
\begin{align*}
    &(1-6) &&\delta_i\delta_j, &&i\neq j,\\
    &(7-18) &&\delta_{ij}\delta_{ih}, &&i\neq j\neq h, \\
    &(19-30)     &&\delta_{i}\delta_{hj}, && i\neq j\neq h, \\
    &(31-42)     &&\delta_{i}(\delta_{ij}-\delta_{ih}), &&i\neq j \neq h, \\
    &(43-48) &&\delta_{ij}(\lau+\delta_{ij}+\delta_{i}+\delta_{j}), &&i\neq j, \\
    &(49-60) &&\delta_i(\lau+\delta_{ij}+\delta_i), &&i \neq j, \\
    &(61-72)      &&\delta_i\delta_{ij}^2, &&i\neq j \\
    &(73-75) &&12(\delta_{ij}^3-\delta_{ik}^3-\delta_{jh}^3+\delta_{hk}^3), &&\{i,j,h,k\}=[4],\\
    &(76-87) &&12(2\delta_{i}^3+6\delta_{i}^2\delta_{ik}+\delta_{ik}^3-\delta_{ik}\delta_{jh}^2-\delta_{ih}\delta_{jk}^2-\delta_{jk}^3-\delta_{ij}\delta_{hk}^2-\delta_{hk}^3), &&\{i,j,h,k\}=[4],\\
    &(88-93) && f_{ij} \text{ (see below)}, &&\{i,j,h,k\}=[4], \\
\end{align*}
\vspace{-1cm}
\begin{align*}
    f_{ij}=&4(\lau^3 + \delta_{i}^3 + \delta_{j}^3 + \delta_{h}^3 + \delta_{k}^3 + \delta_{ij}^3 + \delta_{ih}^3 \\
    &+ 3\delta_{i}^2\delta_{ik} - 2\delta_{ik}^3 - \delta_{ik}\delta_{jh}^2 - 2\delta_{jh}^3 + 3\delta_{j}^2\delta_{jk} - \delta_{ih}\delta_{jk}^2 + \delta_{jk}^3 + 3\delta_{h}^2\delta_{hk} + 3\delta_{k}^2\delta_{hk} - \delta_{ij}\delta_{hk}^2 + \delta_{hk}^3).
\end{align*}
The Hilbert series of the rational Chow ring of $\overline{M}_{1,4}^{ps}$ is equal to $1+11t+18t^2+11t^3+t^4$.
\end{proposition}

\subsection{The Chow ring of $\oM_{1,4}$}
We know from \Cref{cor:picture M14} that $\oM_{1,4}$ is a weighted blow-up of $\oM_{1,4}(1)$ centered at $\Cu$, which is the closed substack of cuspidal curves. As usual, we need the following three ingredients in order to apply the weighted blow-up formula for Chow rings:
\begin{enumerate}
    \item surjectivity and the kernel of the pullback homomorphism $\iota^*\colon A^*(\oM_{1,4}(1))\to A^*(\oM_{0,5})$.
    \item the fundamental class of $\Cu\subset\oM_{1,4}(1)$,
    \item the normal bundle of $\Cu$, and
\end{enumerate}

Recall that the locus of cuspidal curves $\Cu_{1,4}$ is isomorphic to $\oM_{0,5}$ by pinching the fifth section.  In turn, $\oM_{0,5}$ is isomorphic to the blow-up of $\PP^1\times\PP^1$ in three points $(0,0)$, $(1,1)$ and $(\infty,\infty)$. The calculation of its Chow rings follows from this explicit description. We take it directly from \cite[Theorem 1]{Keel}: let $S$ be a subset of the markings set $[5]$ of cardinality between $2$ and $3$, and let $D^S$ be (the class of) the divisor in $\oM_{0,5}$ of curves with a node separating the markings indexed by $S$ from those indexed by $S^c$. We say that two sets of indices $S$ and $T$ are incomparable if $S\nsubseteq T, S^c\nsubseteq T,T\nsubseteq S, T^c\nsubseteq S$. Keel's result is the following:
\begin{theorem}[Keel]
The integral Chow ring of $\oM_{0,5}$ is :
    \[A^*(\oM_{0,5})=\ZZ[D^S]/\left(D^S-D^{S^c},\sum_{i,j\in S,h,k\in S^c}D^S-\sum_{i,h\in S,j,k\in S^c}D^S; D^SD^T\right)_{\substack{h,i,j,k\in[5]\\S,T\text{ incomparable}}}\]
\end{theorem}
We observe that the following relations hold under pullback:
\begin{align}\label{eqn:bdry_restr}
    \iota^*\delta_{hk}=D^{\{i,j\}},&\qquad \{i,j\}\subseteq [4],\\
    \iota^*\delta_i=D^{\{i,5\}},&\qquad i\in [4].
\end{align}
This already shows the surjectivity of $\iota^*$.

To compute the restriction of $\lau$, we first relate it to a $\psi$-class on $\oM_{1,4}(1)$, and then express every term in terms of Keel's basis after pulling them back to $\oM_{0,5}$. First, \Cref{prop:psiVSlambda} identifies the difference between $\lau$ and $\psi_1$ with the class of the divisor where the first marking lies on a rational tail:
\begin{equation}
    \lau=\psi_1-\sum_{\{i,j,k\}=\{2,3,4\}}(\delta_{ij}+\delta_k).
\end{equation}
Observe now that $\iota^*\psi_1=\psi_1$ because the pinching occurs away from the first marking. Finally, we write $\psi_1\in A^1(\oM_{0,5})$ in terms of boundary classes. The comparison of $\psi$ classes under forgetting the last marking gives us:
\begin{equation}
    \psi_1^{\{0,5\}} = f_5^*\psi_1^{\{0,4\}} + D^{\{1,5\}}.
\end{equation}
For the same reason, $\psi_1^{\{0,4\}}=D^{\{1,4\}}$, so we conclude:
\begin{align*}
    \iota^*\lau &=D^{\{1,4\}}+D^{\{2,3\}}+D^{\{1,5\}}-(D^{\{1,2\}}+D^{\{1,3\}}+D^{\{1,4\}}+D^{\{2,5\}}+D^{\{3,5\}}+D^{\{4,5\}})\\
    &= -(D^{\{1,3\}}+D^{\{2,5\}}+D^{\{4,5\}})=-\iota^*(\delta_{24}+\delta_2+\delta_4).
\end{align*}
Putting this together with \Cref{eqn:bdry_restr}, we deduce the following.
\begin{lemma}\label{lm:ker cusp}
The kernel of the pullback homomorphism $\iota^*\colon A^*(\oM_{1,4}(1)) \to A^*(\oM_{0,5})$ is generated by:
    \begin{itemize}
        \item[(1)] $(\delta_{12}+\delta_{34}) - (\delta_{13}+\delta_{24})$,
        \item[(2)] $(\delta_{12}+\delta_{34}) - (\delta_{14}+\delta_{23})$,
        \item[(3-10)] $(\delta_{hk}+\delta_h) - (\delta_{jk}+\delta_j),\qquad j,h,k\in[4]$,
        \item[(11)] $\lau + \delta_{24} + \beta_2 + \beta_4$.
    \end{itemize}
\end{lemma}

Finally, from \Cref{cor:Ell_classes} we know that:
\begin{align*}
    [\Cu]&=24\lau^2\in A^2(\oM_{1,4}(1)),\\
    \Ncal_{\Cu/\oM_{1,4}(1)}&=(\Hcal(1)^{\otimes 4}\oplus \Hcal(1)^{\otimes 6})_{|\Cu}.
\end{align*}
Denoting by $\delta_{\emptyset}$ the fundamental class of the divisor in $\oM_{1,4}$ of curves with an unmarked elliptic tail, applying \eqref{eqn:blow-up formula} gives us the following.
\begin{proposition}
    We have 
    \[A^*(\oM_{1,4})=\ZZ[\lau,\{\delta_{ij}\}_{\{i,j\}\subset [4]},\{\delta_i\}_{i\in [4]},\delta_{\emptyset}]/I_{1,4}'\]
    where $I_{1,4}'$ is generated by the relations
    \begin{align*}
    & (1-93) && \text{as in \Cref{prop:chow M14(1)},} && \\
    &(94) &&\delta_{\emptyset}((\delta_{12}+\delta_{34}) - (\delta_{13}+\delta_{24})), && \\
    &(95) &&\delta_{\emptyset}((\delta_{12}+\delta_{34}) - (\delta_{14}+\delta_{23})), && \\
    &(96-103) &&\delta_{\emptyset}((\delta_{hk}+\delta_h) - (\delta_{jk}+\delta_j)), && j,h,k\in[4],\\
    &(104) &&\delta_{\emptyset}(\lau + \delta_{24} + \delta_2 + \delta_4), && \\
    &(105) &&24(\lau-\delta_{\emptyset})^2. &&
\end{align*}
\end{proposition}
To conclude, because of \Cref{prop:lambdaVSlambda} we have $\lambda=\lau-\delta_{\emptyset}$, where $\lambda$ is the first Chern class of the Hodge line bundle on $\oM_{1,4}$. Applying this change of variable, we finally get to the following presentation for $A^*(\oM_{1,4})$. As above, we do not know whether it is minimal.
\begin{theorem}\label{thm:chow M14}
    We have
    \[A^*(\oM_{1,4})=\ZZ[\lambda,\{\delta_{ij}\}_{\{i,j\}\subset [4]},\{\delta_i\}_{i\in [4]},\delta_{\emptyset}]/I_{1,4}\]
    where $I_{1,4}$ is generated by the following relations:
    \begin{align*}
    &(1-6) &&\delta_i\delta_j, &&i\neq j,\\
    &(7-18) &&\delta_{ij}\delta_{ih}, &&j\neq h, \\
    &(19-30)     &&\delta_{i}\delta_{hj}, && h,j \neq i \\
    &(31-42)     &&\delta_{i}(\delta_{ij}-\delta_{ih}), &&i\neq j \neq h, \\
    &(43-48) &&\delta_{ij}(\lambda + \delta_{\emptyset} +\delta_{ij}+\delta_{i}+\delta_{j}), &&i\neq j, \\
    &(49-60) &&\delta_i(\lambda + \delta_{\emptyset}+\delta_{ij}+\delta_i), &&j,h \neq i, \\
    &(61-72)      &&\delta_i\delta_{ij}^2, &&i\neq j \\
    &(73-75)  &&12(\delta_{ij}^3-\delta_{ik}^3-\delta_{jh}^3+\delta_{hk}^3), &&\{i,j,h,k\}=[4],\\
    &(76-87) &&12(2\delta_{i}^3+6\delta_{i}^2\delta_{ik}+\delta_{ik}^3-\delta_{ik}\delta_{jh}^2-\delta_{ih}\delta_{jk}^2-\delta_{jk}^3-\delta_{ij}\delta_{hk}^2-\delta_{hk}^3), &&\{i,j,h,k\}=[4],\\
    &(88-93) && f_{ij}\text{ (see below)}, &&\{i,j,h,k\}=[4]. \\
    &(94) &&\delta_{\emptyset}((\delta_{12}+\delta_{34}) - (\delta_{13}+\delta_{24})), && \\
    &(95) &&\delta_{\emptyset}((\delta_{12}+\delta_{34}) - (\delta_{14}+\delta_{23})), && \\
    &(96-103) &&\delta_{\emptyset}((\delta_{hk}+\delta_h) - (\delta_{jk}+\delta_j)), && j,h,k\in[4],\\
    &(104) &&\delta_{\emptyset}(\lambda + \delta_{\emptyset} + \delta_{24} + \delta_2 + \delta_4), && \\
    &(105) &&24\lambda^2, &&
\end{align*}
\vspace{-1cm}
\begin{align*}
    f_{ij}=&4((\lambda+\delta_{\emptyset})^3 + \delta_{i}^3 + \delta_{j}^3 + \delta_{h}^3 + \delta_{k}^3 + \delta_{ij}^3 + \delta_{ih}^3 \\
    &+ 3\delta_{i}^2\delta_{ik} - 2\delta_{ik}^3 - \delta_{ik}\delta_{jh}^2 - 2\delta_{jh}^3 + 3\delta_{j}^2\delta_{jk} - \delta_{ih}\delta_{jk}^2 + \delta_{jk}^3 + 3\delta_{h}^2\delta_{hk} + 3\delta_{k}^2\delta_{hk} - \delta_{ij}\delta_{hk}^2 + \delta_{hk}^3).
\end{align*}
The Hilbert series of the rational Chow ring of $\overline{M}_{1,4}$ is equal to $1+12t+23t^2+12t^3+t^4$.
\end{theorem}

\appendix
\section{Logarithmic resolution of forgetful maps}\label{appendix}
We come back to the problem of \Cref{sec:tilda} in greater generality and in a different language.

Recall that, contrary to the nodal case, the forgetful map $\oM_{1,n}(m)\dashrightarrow\oM_{1,n-1}(m)$ is only rational for $m\geq2$; in other words, $\oM_{1,n}(m)$ is only birational to the universal curve over $\oM_{1,n-1}(m)$. Basically, the reason is that elliptic $m$-fold points - as well as their higher genus generalisations - have non-trivial moduli of attaching data, so these singularities are not determined by their pointed normalisation.

Suppose that the $n$-th marking lies on the core, and that the latter has level $m+1$. Then forgetting $p_n$ makes the curve $m$-unstable. The solution is to contract the core after sprouting the markings lying directly on it; this will increase the level, since there must be at least one rational tail containing at least two special points.
The contraction is not well-defined as soon as there are at least two rational tails. In order to resolve the indeterminacy, it is enough to blow up the loci $\EB_S(m)$, where $S$ is a partition of $[n]$ into $m+1$ parts, of which at least two have $\lvert S_i\rvert\geq2$, and one is the singleton $\{n\}$; see \cite[\S 2.1]{SmythIII}.

We now want to forget $n-m$ markings at once. Smyth's procedure would produce a zig-zag of spaces. With a touch of logarithmic geometry we may instead resolve the indeterminacy simultaneously by adapting the strategy of \cite{RSPW}, as we now explain.

To every $n$-marked $m$-stable elliptic curve $C$ we can associate an \ul{$m$-marked} genus-weighted graph $\Gamma_m$ as follows: $\Gamma_m$ has one vertex (of genus $1$) corresponding to the core, and one vertex (of genus $0$) for every irreducible component not belonging to the core; edges correspond to nodes, except for those internal to the core (in case the latter is a cycle of $\PP^1$s); moreover, $\Gamma_m$ has $m$ legs corresponding to \emph{the first $m$ markings} of $C$, each attached to the vertex corresponding to the subcurve containing it. Finally, vertices of genus $0$ and valence $1$ or $2$ are removed as follows: bridges are forgotten, while for tails we perform an edge contraction of the unique adjacent edge.

We now define $\oM_{1,n}(m)^\dagger$ as the moduli functor over log schemes representing \emph{aligned} $n$-marked $m$-stable curves of arithmetic genus one: these curves have at worst $m$-fold elliptic singularities; they are endowed with the trivial log structure on the complement of the rational tails (the "interior" of the core), and with the log smooth log structure away from the singular points of the core, so that for every geometric point $s$ of the base $S$ the dual graph $\Gamma_m(C_s)$ is naturally metrised in the base characteristic monoid $\overline{M}_{S,s}$. Finally, we require that the tropical distances from the rational tails to the core are totally ordered in the base monoid (alignment).

\begin{proposition}
 Let $m\leq 5$. Then $\oM_{1,n}(m)^\dagger$ is represented by a \emph{smooth} Deligne--Mumford stack which is log smooth with a locally free logarithmic structure. It is a log blow-up of $\oM_{1,n}(m)$.
\end{proposition}
\begin{proof}
    Analogous to \cite[Proposition 3.3.4]{RSPW}. Consider the forgetful map $\oM_{1,n}(m)\to \Ucal_{1,m}$, where the latter denotes the stack of $m$-marked Gorenstein curves of genus one. As usual, we can replace $\Ucal_{1,m}$ with a finite type substack containing the image of $\oM_{1,n}(m)$ (in particular, we can restrict to curves with at worst elliptic $m$-fold points); we abuse notation by assuming this done. Now, endow $\Ucal_{1,m}$ with the toroidal log structure induced by the divisor of rational tails. Aligning corresponds to a logarithmic blow-up of the latter, as explained in \cite[\S 3.4]{RSPW}, which is then pulled back to $\oM_{1,n}(m)$.
\end{proof}

The Picard group of $\oM_{1,n}(m)$ is freely generated by $\lambda$ and boundary divisors $\tau_S$, $S\subseteq[n]$ such that $2\leq\lvert S\rvert\leq n-m$ denoting the markings on the rational tail \cite[Proposition 3.2]{SmythII} (recall that $\lambda$ and $\delta_{\text{irr}}$ are multiple of one another). On $\oM_{1,n}(m)^\dagger$ we have the exceptional divisors as well \cite[Proposition 3.2]{SmythIII}. They can be interpreted modularly and described in tropical terms. In fact, the Picard group of $\oM_{1,n}(m)^\dagger$  is generated by those boundary strata that depend on a single tropical parameter. These strata $\tau_{\Scal}$ are now indexed by partitions $\Scal$ of the marking set $[n]$: two markings belong to the same part if there is a node separating them from the core. In particular, $\Scal$ has at least $m+1$ (level condition) and at most $n-1$ parts.

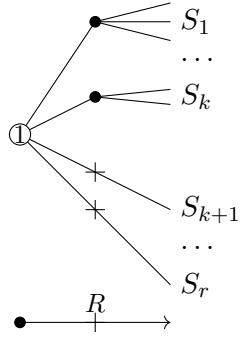
\begin{figure}[hbt]
 \begin{tikzpicture}
  \tikzset{cross/.style={cross out, draw=black, fill=none, minimum size=2*(#1-\pgflinewidth), inner sep=0pt, outer sep=0pt}, cross/.default={3pt}}

   \foreach \x in {(1,1.5),(1,.5),(0,-2.5)}
  \draw[fill=black] \x circle (2pt);

  \foreach \x in {(1,-1),(1,-.5),(1,-2.5)}
  \draw \x node[cross,rotate=45] {};

  \draw (0,0) -- (1,1.5) -- (2,1.75) (1,1.5) -- (2,1.5) node[right]{$S_1$} (1,1.5) -- (2,1.25)
  (2,1) node[right]{$\ldots$}
  (0,0) -- (1,.5) -- (2,.6) (1,.5)--(2,.4) (2,.5) node[right]{$S_k$}
  (0,0) -- (1,-.5)--(2,-1) node[right]{$S_{k+1}$}
  (2,-1.5) node[right]{$\ldots$}
  (0,0)--(1,-1)--(2,-2) node[right]{$S_r$};

  \draw[fill=white] (0,0) circle (4pt) node {\small $1$};
  \draw[->] (0,-2.5)--(2,-2.5);
  \draw (1,-2.5) node[above]{\small $R$};
 \end{tikzpicture}
\caption{The generic dual graph in the divisor $\tau_{\Scal}$: there is at least one stable rational tail, and all rational tails are at the same distance $R$ from the core ($r>m$).}
\label{fig:tropdiv}
\end{figure}

Let us denote by $b\colon \oM_{1,n}(m)^\dagger\to \oM_{1,n}(m)$ the blow-up map. Then, $\tau_{\Scal}$ is an exceptional divisor of $b$ precisely when there is more than one part of $\Scal$ of cardinality $\geq2$ (the number of such parts is the codimension of the corresponding stratum in $\oM_{1,n}(m)$); at the generic point of $\tau_{\Scal}$, all the stable rational tails are equidistant from the core. See Figure \ref{fig:tropdiv}.

\begin{proposition}\label{prop:forgetful}
    There is a forgetful map $f_{n,m}\colon \oM_{1,n}(m)^\dagger\to \widetilde{\mathcal{M}}_{1,m}$ such that:
    \begin{equation}\label{eqn:forgetful_pullback}
        f_{n,m}^*(\widetilde{\Hcal}(m))=\Hcal(m)\otimes\OO_{\oM_{1,n}(m)^\dagger}\left(\sum_{\Scal\vdash [n],\lvert \Scal\rvert\geq m+1,[m]\to \Scal\text{ not injective}}\tau_{\Scal}\right).
    \end{equation}
    Moreover, the preimage of the $m$-fold point $\Zcal$ is the union of the loci of $m$-fold points such that every $p_i,\ i=1,\ldots,m,$ cleaves to a different branch $\mathbf{Ell}_{\{\{p_1,\ldots\},\ldots,\{p_m,\ldots\}\}}(m)$. In a neighbourhood of this locus, $b$ is an isomorphism, and $f_{n,m}$ is flat; in particular, the normal bundle can be read off from \Cref{lm:class of Z}.
\end{proposition}
\begin{proof}
    The construction of the forgetful map follows along the lines of \cite[Theorem 3.7.1]{RSPW}: consider circles or radius $R$ around the core; let $R_m$ be the minimal radius such that the inner valence of the corresponding circle is $\leq m$ (this determines the singularity of the contraction), and the outer valence is $\geq m$ (this determines the core level of the contraction). There is a log modification $\widetilde C$ of the universal curve $C$ over $\oM_{1,n}(m)^\dagger$ corresponding to the subdivision of every edge (finite or infinite) at distance $R_m$ from the core, cf. \cite[Proposition 3.6.1]{RSPW}; geometrically, this corresponds to sprouting markings and blowing up nodes as necessary to make sure that, when marked with the strict transforms of $p_1,\ldots,p_m$, the markings may lie on the core only if the level of the latter is at least $m$.

    Let $\mu$ denote the conewise-linear (CL) function on $\Gamma_m$ representing the minimum between $R_m$ and the distance from the core: the fundamental exact sequence of log geometry associates to $\mu$ a line bundle $\OO_{\widetilde C}(\mu)$, which can be interpreted fibrewise according to the dictionary of \cite[Proposition 2.4.1]{RSPW}. In particular, the line bundle $\Lcal=\omega_{\widetilde C}(p_1+\ldots+p_m-\mu)$ is relatively semiample, trivial on the disc of radius $R_m$, and ample outside of it, so that we obtain a contraction:
    \bcd
    \widetilde C\ar[rr,"\tau"]\ar[dr] & & \overline C\ar[dl] \\
    & \oM_{1,n}(m)^\dagger &
    \ecd
    where $\overline C=\on{Proj}_{\Mcal^\dagger}(\bigoplus_{k\geq0}\tilde\pi_*\Lcal^{\otimes k})$ is almost $(m-1)$-stable when pointed with $p_1,\ldots,p_m$, cf. \cite[Proposition 3.7.6.1]{RSPW}. Roughly speaking, $\tau$ performs the following operations:
    \begin{itemize}
        \item forgets all but the first $m$ markings,
        \item contracts all the rational tails and bridges that have thus been made unstable,
        \item contracts a balanced subcurve of genus one in such a way to produce an elliptic $l$-fold point, with $l\leq m$, and level at least $m$.
    \end{itemize}
    This contraction induces the morphism $f_{n,m}\colon \oM_{1,n}(m)^\dagger\to \widetilde{\mathcal{M}}_{1,m}$.

    Proposition \ref{prop:lambdaVSlambda} implies that the Hodge bundle is unaffected by the contraction of rational components, whereas it changes by the divisorial loci over which the core is contracted. This happens precisely when the level of the core of $C$, marked by $p_1,\ldots,p_m$ only, is strictly lower than $m$, i.e. when at least two of the first $m$ markings belong to the same rational tail. This is exactly the condition that the map $[m]\to S$ associating to $i$ the part to which $p_i$ belongs is not injective. 

    Finally, we notice that, as soon as two among the first $m$ markings lie on the same rational tail, the above procedure produces an elliptic $l$-fold point with $l<m$ in $\overline C$; hence, $\overline C$ has an elliptic $m$-fold point if and only if $C$ has one, such that each $p_i$ cleaves (possibly through a rational chain) to a distinct branch. This locus is contained in the complement of the divisor where at least two markings among the first $m$ are separated from the core by a (the same) separating node, and in particular it is disjoint from the blow-up centre.

    The flatness of $F$ can be proved by induction on $n$: the base case is $n=m+1$. In this case, the morphism $F\colon \oM_{1,m+1}(m)\to\widetilde{\mathcal{M}}_{1,m}$ is defined without need of blowing up, because the curve always coincides with its core. In particular, no core contraction ever occurs, so that $F$ identifies $\oM_{1,m+1}(m)$ with the universal curve over $\widetilde{\mathcal{M}}_{1,m}$ minus the universal sections and the $l$-fold elliptic points; being the composition of an open embedding and a flat family, $F$ is therefore flat. The inductive step follows along the same lines.
\end{proof}

\bibliographystyle{halpha-abbrv}
\bibliography{bibliography.bib}

\subsection*{Funding} During the preparation of this work, L.B. was partially supported by the ERC Advanced Grant SYZYGY of the European Research Council (ERC) under the European Union Horizon 2020 research and innovation program (grant agreement No. 834172) and by the European Union - NextGenerationEU under the National Recovery and Resilience Plan (PNRR) - Mission 4 Education and research - Component 2 From research to business - Investment 1.1 Notice Prin 2022 - DD N. 104 del 2/2/2022, from title "Symplectic varieties: their interplay with Fano manifolds and derived categories", proposal code 2022PEKYBJ – CUP J53D23003840006. L.B. is a member of INdAM group GNSAGA.


\end{document}